\documentclass[11pt]{article}
\usepackage{amssymb,latexsym}
\usepackage{amsmath}

\title{ \large {\bf Finitely generated algebras with involution
and their identities.}}

\author{ Irina Sviridova\thanks{Supported by CNPq, DPP/UnB and by CNPq-FAPDF PRONEX
grant 2009/00091-0 (193.000.580/2009); e-mail \texttt{I.Sviridova@mat.unb.br}}
\\
\\
Departamento de Matem\'atica,\\
Universidade de Bras\'\i lia,\\
70910-900 Bras\'\i lia, DF, Brazil }

\date{February 11, 2013.}

\newtheorem{theorem1}{Theorem}
\newtheorem{theorem}{Theorem}
\newtheorem{definition}{Definition}
\newtheorem{lemma}{Lemma}

\newtheorem{remark}{Remark}

\begin{document}
\maketitle

\begin{abstract}
We consider associative algebras with involution over a field of
characteristic zero. In this case, we prove that for any
finitely generated associative algebra with involution there
exists a finite dimensional algebra with involution which
satisfies the same identities with involution. This is an analogue
and an extension of the theorem of  A.Kemer for ordinary
identities \cite{Kem1}. The similar results were proved earlier
by the author for identities graded by a finite abelian group
\cite{Svi}, and by E.Aljadeff, and
A.Kanel-Belov \cite{AB} for identities graded by any finite group.

\bigskip


\end{abstract}

\section*{Introduction}

The theory of polynomial identities of algebras is an important area of the modern algebra.
In the last years various generalizations of polynomial identities become
popular. The purpose of the present paper is to study identities
with involution. We prove that finitely generated algebras with involution
are equivalent to finite dimensional algebras in terms of $*$-identities,
the field being of characteristic zero.

Note that almost all notions mentioned here can be applied not
only for associative algebras but for other classes of algebras.
Nevertheless throughout the paper we
consider only associative algebras over a field of
characteristic zero. Further they will simply be called algebras.

Let $F$ be a field of characteristic zero. We consider
associative algebras over $F$. An anti-automorphism $*$ of
the second order of an $F$-algebra $A$ is called an {\it
involution}. If we fix an involution $*$ of an associative
$F$-algebra $A$ then the pair $(A,*)$ is called an {\it associative
algebra with involution} (or {\it associative $*$-algebra}). Note that an
algebra with involution can be considered as an algebra with the
supplementary unary operation $*$ satisfying identities
\begin{eqnarray*}
&&(\alpha a+ \beta b)^*=\alpha a^* + \beta b^*, \\ &&(a \cdot
b)^*=b^* \cdot a^*, \qquad (a^*)^*=a
\end{eqnarray*}
for any $a, b \in A,$ and $\alpha, \beta \in F.$

Observe that any $*$-algebra can be decomposed into
the sum of symmetric and skew-symmetric parts. An element $a \in
A$ is called symmetric if $a^*=a,$ and
skew-symmetric if $a^*=-a.$ So, $a+a^*$ is symmetric and
$a-a^*$ skew-symmetric for any $a \in A.$ Thus, we have
$A=A^{+} \oplus A^{-},$ where $A^{+}$ is the subspace
formed by all symmetric elements ({\it symmetric part}), and
$A^{+}$ is the subspace of all skew-symmetric elements of $A$
({\it skew-symmetric part}). We call symmetric and
skew-symmetric elements of a $*$-algebra {\it $*$-homogeneous}
elements.

We use standard notations $a \circ b =a b + b a,$ and $[a,b]=a b - b a.$
It is well known that the symmetric part $A^{+}$ of
a $*$-algebra $A$ with the operation $\circ$ is a Jordan algebra, and
the skew-symmetric part $A^{-}$ with the operation $[,]$ is a Lie algebra.

A homomorphism of algebras with involution is a
$*$-homomorphism if it commutes with the involution. An ideal $I
\unlhd A$ of a $*$-algebra $A$ is $*$-ideal if it is invariant
under the involution. For algebras with involution we
consider only $*$-ideals, and $*$-homomorphisms. In this case the
quotient algebra $A/I$ is also an algebra with involution.

We denote by $A_1 \times \dots \times A_{\rho}$ the direct product
of algebras $A_1, \dots, A_{\rho},$ and by $A_1 \oplus \dots
\oplus A_{\rho} \subseteq A$ the direct sum of subspaces $A_i$ of
an algebra $A.$ It is clear that a direct product of algebras
with involution is also an algebra with involution. Throughout the
paper we denote by $J(A)$ the Jacobson radical of $A.$ By default,
all bases and dimensions of spaces and algebras are considered
over the base field $F.$ We consider the lexicographical order on
$n$-tuples of numbers.

The notion of identity with involution (a $*$-identity) is a formal extension of
the notion of ordinary polynomial identity.
We refer the reader to the textbooks \cite{Dren},
\cite{DrenForm}, \cite{GZbook}, and to \cite{BelRow}, \cite{Kem1} concerning basic
definitions, facts and properties of polynomial identities.

A brief introduction to identities with involution
will be given in Section 1 (see also \cite{GZbook}).
Section 2 contains an information about basic properties of
the main objects and constructions. Here we also introduce the parameters $\mathrm{par}_*(A)$
of a finite dimensional algebra $A$, and the Kemer index $\mathrm{ind}_*(\Gamma)$
of an ideal $\Gamma$ of identities with involution of a finitely generated algebra
with involution.  Section 3 is devoted to
finite dimensional $*$-algebras, namely, to connection between their structural
parameters $\mathrm{par}_*,$ and the indices of their ideals
of $*$-identities $\mathrm{ind}_*.$
Section 4 is devoted to representable algebras, and to the technique
of approximation of finitely generated algebras with involution
by finite dimensional $*$-algebras. The concept of form identities
with involution is also considered in this section. Section 5 contains
the proof of the following main result.

\begin{theorem1} \label{*PI}
Let $F$ be a field of zero characteristic. Then
a non-zero $*T$-ideal of \mbox{$*$-identities} of a finitely generated
associative $F$-algebra with involution coincides
with the $*T$-ideal of \mbox{$*$-identities} of some finite dimensional
associative $F$-algebra with involution.
\end{theorem1}

The methods and techniques are an adaptation of the methods
used by the author in the proof of the corresponding result for
graded identities of finitely generated algebras
graded by a finite abelian group \cite{Svi}, see also \cite{AB}
for a more general result. In both cases the
origin of methods is the techniques developed
by Alexander Kemer to obtain
the classification of varieties of non-graded and
$\mathbb{Z}_2$-graded associative algebras and the positive
solution of the Specht problem (\cite{Kem1}).
The useful interpretation of his methods can be found in
\cite{BelRow}.

In general, we follow the structure of the proof \cite{Svi}
with necessary modifications for the case of involution.
Note that for several technical statements the proof is the literal repetition
of the corresponding proof for the graded identities. In such cases we
refer the reader to \cite{Svi}.
The origin of the principal notations and definitions is \cite{BelRow} (see also
\cite{AB}). More details concerning the basic constructions and
definitions also can be found there.

\section{Free algebra with involution.}

Consider a countable set of letters $X=\{ x_{i} | i \in \mathbb{N} \}$, and denote
$X^* = \{ x_{i}, x_{i}^* | x_i \in X \}.$ We define $*$-action
on monomials in $X^*$ by
\begin{eqnarray*}
&&(a_{i_1} \cdots a_{i_n})^*=a_{i_n}^* \cdots a_{i_1}^*, \quad a_{j} \in X^*, \quad
\mbox{ where } \quad (x_j^*)^*=x_j \quad \mbox{for all } \ x_j \in X.
\end{eqnarray*}
This action is naturally extended to the involution
on the free associative algebra $\mathfrak{F}=F\langle X^{*}
\rangle$ generated by the set $X^*.$ The algebra $\mathfrak{F}$
with this involution is called the free associative algebra with
involution. Then an element of the free algebra with involution
has a form $f(x_1,\dots,x_n)=\sum_{(i), (\theta)} \alpha_{(i),
(\theta)} x_{i_1}^{\theta_1} \cdots x_{i_m}^{\theta_m},$ where
$x_i \in X,$ $\alpha_{(i), (\theta)} \in F,$ \ $\theta_i \in \{ 0,
1\},$ and $x_i^{\theta}=x_i$ \ if $\theta=0,$ \
$x_i^{\theta}=x_i^*$ \ if $\theta=1.$

Taking into account the decomposition of $*$-algebras into
symmetric and skew-symmetric parts we can consider another
construction of the free associative algebra with involution. Consider
two countable sets $Y = \{ y_{i} | i \in \mathbb{N} \},$ \ $Z = \{ z_{i}
| i \in \mathbb{N} \}.$  The action $*$ is defined on monomials in $Y
\cup Z$ by equalities
\begin{eqnarray*}
&&w^*=(a_{i_1} \cdots a_{i_n})^*=a_{i_n}^* \cdots a_{i_1}^*
=(-1)^{\delta(w)} a_{i_n} \cdots a_{i_1}, \quad \mbox{where}\\
&&y_j^*=y_j, \ \ z_j^*=-z_j, \quad a_j \in Y \cup Z.
\end{eqnarray*}
Here the sign is determined by the parity of the number
$\delta(w)$ of variables from the set $Z$ in the monomial $w.$
The linear extension of this action is an involution on the free
associative algebra $F\langle Y, Z \rangle$ generated by the set
$Y \cup Z.$ The equalities
\begin{eqnarray} \label{free}
&&y_i=\frac{x_i+x_i^*}{2}, \ \  z_i=\frac{x_i-x_i^*}{2}; \nonumber
\\ &&x_i=y_i+z_i, \ \ x_i^*=y_i-z_i
\end{eqnarray}
induce an isomorphism of algebras with involution $F\langle
X^*\rangle$ and $F\langle Y, Z \rangle.$ We use both
notations for the free associative algebra with involution
$\mathfrak{F}.$ We call a set of variables $S \subseteq Y$ or
$S \subseteq Z$ (namely, if $S$ does not contain variables from $Y$ and $Z$
simultaneously) {\it $*$-homogeneous}.

Note that degrees of polynomials in the first case are defined by
conditions $\deg_{x_i} x_i=\deg_{x_i} x_i^*=1$ ($i \in
\mathbb{N}$).

Let $f=f(x_{1}, \dots, x_{n}) \in F\langle X^{*} \rangle$ be a
non-trivial $*$-polynomial. We say that a $*$-algebra $A$
satisfies {\it the identity with involution } (or {\it $*$-identity})
$f(x_{1},\dots x_{n})=0$ if $f(a_{1}, \dots a_{n})=0$ for
any $a_{i} \in A.$ If a polynomial $f=f(y_1,\dots,y_n,z_1,\dots,z_m)
\in F\langle Y, Z \rangle$ is written in terms of the symmetric and
skew-symmetric variables then $f=0$ is $*$-identity of $A$ iff
$f(a_1,\dots,a_n,b_1,\dots,b_m)=0$ for all $a_i \in A^{+},$ and
$b_i \in A^{-}.$

Let $\mathrm{Id}^{*}(A) \unlhd F\langle X^{*} \rangle$ be the
ideal of all identities with involution of $A.$ Similar to the
case of ordinary identities any ideal of identities with
involution is a two-sided \mbox{$*$-ideal} of the free algebra with
involution $F\langle X^{*} \rangle$ invariant under its
\mbox{$*$-endomorphisms}. We call such ideals {\it $*$T-ideals}. Conversely, any
$*$T-ideal of $F\langle X^{*}\rangle$ is the ideal of identities
with involution of some $F$-algebra with involution.

Given a monomial $w=a_1 a_2 \cdots a_n$ in $X^{*}$ (in case
$a_i \in X^{*}$) or in $Y \cup Z$ (if $a_i \in Y \cup Z$) let us
denote by $\widetilde{w}=a_n \cdots a_2 a_1$ the monomial rewritten in
the reverse order. Given a polynomial $f=\sum_w \alpha_w w,$ we
denote $\tilde{f}=\sum_w \alpha_w \widetilde{w}$ (where $w$-s are
monomials, $\alpha_w \in F$). Note that this operation is also
an involution of $\mathfrak{F}.$ It is clear that for any $f \in
\mathfrak{F}$ we have
$f(x_1,\dots,x_n)^*=\tilde{f}(x_1^*,\dots,x_n^*).$ In particular,
polynomials $f$, and $\tilde{f}$ always belong to the same
$*$T-ideal.

For a set $S \subseteq F\langle X^{*} \rangle$ of $*$-polynomials
denote by $*T[S]\unlhd F\langle X^{*} \rangle$ the
$*$T-ideal generated by $S.$ Then $*T[S]$ contains exactly all
identities with involution which follow from polynomials of the
set $S.$ We say that two algebras with involution $A$ and $B$ are
$*$PI-equivalent $A \sim_{*} B$ if
$\mathrm{Id}^*(A)=\mathrm{Id}^*(B).$ We also write $f=g \
(\mathrm{mod } \ \Gamma)$ for a $*$T-ideal $\Gamma$ and
$*$-polynomials $f, g \in F\langle X^{*} \rangle$ if $f-g \in
\Gamma.$

It is clear that the ordinary free associative algebra $F \langle X
\rangle$ (without involution) can be considered as a subalgebra of
$F \langle X^{*} \rangle.$ Particularly, an ordinary polynomial
identity (without involution) can be considered as an identity with
involution. So, let $A$ be a $*$-algebra, $\mathrm{Id}^{*}(A)$ ideal of
$*$-identities, and $\mathrm{Id}(A)$ ideal of ordinary identities of $A,$
then we have $\mathrm{Id}(A) \subseteq \mathrm{Id}^{*}(A).$
Moreover, by Amitsur's theorem (\cite{A3}, \cite{A4}, see also \cite{GZbook})
any $*$-algebra satisfying a non-trivial $*$-identity also has a
non-trivial ordinary identity.

Note that a finitely generated algebra with involution also
has a finite generating set consisting
of symmetric and skew-symmetric elements. Such set is
called a \mbox{$*$-homogeneous} generating set.

Given a finitely generated $*$-algebra $A,$ we denote by
$\mathrm{rk}(A)$ the least possible number of generators.
Denote by $\mathrm{rkhs}(S)$ the greater of the numbers
of symmetric and skew-symmetric elements of a \mbox{$*$-homogeneous}
generating set $S$ of $A.$ Then $\mathrm{rkh}(A)$ is
the least possible $\mathrm{rkhs}(S)$
for the algebra $A.$

Let $X_\nu^{*} = \{ x_{i}, x_{i}^{*} | 1 \le i \le \nu \}$ be a
finite set and $F\langle X_\nu^{*} \rangle$ the
free associative algebra with involution of rank $\nu,$ \ $\nu
\in \mathbb{N}.$ It is isomorphic to the algebra
$F\langle Y_\nu, Z_{\nu} \rangle$ in symmetric
$Y_\nu = \{ y_{1},\dots,y_{\nu} \},$ and skew-symmetric variables
$Z_\nu = \{ z_{1},\dots,z_{\nu} \}.$
Let $A$ be a finitely generated algebra with involution, then
for any $\nu \ge \mathrm{rk}(A)$
the algebra with involution $U_\nu=F\langle X_\nu^{*}
\rangle/(\mathrm{Id}^*(A) \cap F\langle X_\nu^{*} \rangle)$ is
the relatively free algebra with involution of the rank $\nu,$ and
$\mathrm{Id}^{*}(A)=\mathrm{Id}^{*}(U_\nu).$

Given a $*$T-ideal $\Gamma_1 \subseteq F\langle X^{*} \rangle$ and
a $*$-algebra $B,$ denote by $$\Gamma_1(B)=\{ f(b_1,\dots,b_n) |
f \in \Gamma_1, b_i \in B \} \unlhd B$$ the verbal ideal of $B$
corresponding to $\Gamma_1.$

Similarly to the case of ordinary \cite{Kem1},
and graded identities \cite{Svi},  we have

\begin{remark} \label{freefg}

\begin{enumerate}
\item Let $\Gamma$ be a $*$T-ideal.
A $*$-polynomial $f(x_1,\dots,x_n)$ is an identity with involution
of the relatively free algebra $F\langle X_\nu^{*} \rangle/(\Gamma
\cap F\langle X_\nu^{*} \rangle)$ if and only if
$f(h_1,\dots,h_n) \in \Gamma$ for any $*$-polynomials
$h_1,\dots,h_n \in F\langle X_\nu^{*} \rangle.$
\item Let \ $\Gamma_2=\mathrm{Id}^{*}(A)$ be the ideal of
identities with involution of a finitely generated $*$-algebra
$A$. Then the condition $\Gamma_1(F\langle X_\nu^{*} \rangle)
\subseteq \Gamma_2$ implies $\Gamma_1 \subseteq \Gamma_2$ for any
$\nu \ge \mathrm{rk}(A),$ and the condition
$\Gamma_1(F\langle Y_\nu, Z_\nu \rangle)
\subseteq \Gamma_2$ implies $\Gamma_1 \subseteq \Gamma_2$ for any
$\nu \ge \mathrm{rkh}(A).$
\end{enumerate}
\end{remark}
Note that this remark is true for algebras of any signature.

The notion of degree of a $*$-identity is similar to the case of ordinary
identities (see \cite{Dren, Kem1, GZbook}). In particular, we can
consider multilinear identities.

Let $P_n^*=\mathrm{Span}_F\{ x_{\sigma(1)}^{\theta_1} \cdots
x_{\sigma(n)}^{\theta_n} | \sigma \in \mathrm{S}_n, \theta_i \in
\{ 0, 1\} \}$ be a vector space of all multilinear $*$-polynomials
of degree $n.$
A multilinear $*$-polynomial of degree $n$ has the form
$f(x_{1},\dots,x_{n})= \sum_{\sigma \in \mathrm{S}_n, \ (\theta)
\in \mathbb{Z}_2^{n} } \alpha_{\sigma, (\theta)} x_{\sigma(1)}^
{\theta_{1}} \cdots x_{\sigma(n)}^{\theta_{n}} \in P_n^*.$ Note
that in terms of symmetric and skew-symmetric variables the
space of multilinear polynomials with involution of degree $n$ has
the form $P_n=\mathrm{Span}_F\{ x_{\sigma(1)} \cdots x_{\sigma(n)}
| \sigma \in \mathrm{S}_n, x_i \in Y \cup Z \}$. It is clear that
$P_n$ is the direct sum of subspaces of multihomogeneous and
multilinear polynomials depending on the same symmetric and
skew-symmetric variables. So, a multilinear polynomial
of degree $n$ has the form
$f(y_{1},\dots,y_{k},z_{1},\dots,z_{n-k})= \sum_{\sigma \in
\mathrm{S}_n} \alpha_{\sigma} x_{\sigma(1)} \cdots x_{\sigma(n)}
\in P_{k,n-k},$ where $x_i=y_i$ for $i=1,\dots,k,$ and
$x_i=z_{i-k}$ if $i=k+1,\dots,n$ ($0 \le k \le n$).

Then for an algebra with involution $A$ and for the $*$T-ideal of
its $*$-identities $\Gamma$ we have that $\Gamma_{k,n-k}=P_{k,n-k}
\cap \Gamma$ and $P_{k,n-k}(\Gamma)=P_{k,n-k}(A)=P_{k,n-k}/
(P_{k,n-k} \cap \Gamma)$ are $(FS_{k} \otimes FS_{n-k})$-modules.
Here $S_{k}$ and $S_{n-k}$ act on symmetric and
skew-symmetric variables independently (see, e.g., \cite{GM}).
Moreover, we can assume that the elements of $\Gamma_{k,n-k}$ and
$P_{k,n-k}(\Gamma)$ depend on variables
$\{y_{1},\dots,y_{k}\},$ and $\{ z_{1},\dots,z_{n-k} \}.$ Then the
subspaces of multilinear $*$-identities of $A$ \ $\Gamma_n=P_n
\cap \Gamma,$ and of multilinear parts of relatively free algebra
$P_{n}(\Gamma)=P_{n}(A)=P_{n}/ (P_{n} \cap \Gamma)$ are the direct
sums of ${n}\choose{k}$ copies of $\Gamma_{k,n-k},$ and
$P_{k,n-k}(\Gamma)$ respectively.

It is well known that
in the case of zero characteristic any system of identities
(ordinary or with involution) is equivalent to a system of multilinear
identities. Thus, in the case of zero characteristic it is enough
to consider only multilinear identities.

Then we obtain the next lemma.
\begin{lemma} \label{lemma1}
The ideal of $*$-identities of any finitely generated associative
algebra with involution contains the ideal of $*$-identities of
some finite dimensional algebra with involution.
\end{lemma}
\noindent {\bf Proof.} Let $A$ be a finitely generated associative
algebra with involution, and $\Gamma=Id^{*}(A) \ne 0$ the ideal of
its $*$-identities. From \cite{A3}, \cite{A4} we have that $A$ is
\mbox{PI-algebra}, thus the ideal of its ordinary identities is
non-trivial $Id(A) \ne 0.$ Then from \cite{Kem1} it follows that
there exists a finite dimensional algebra $C$ that has the same
ideal of ordinary identities $Id(A)=Id(C).$ Let us consider the
direct product $B=C \times C^{op},$ where $C^{op}$ is the opposite
algebra. The algebra $B$ has the natural exchange
involution defined by $(a,b)^*=(b,a),$ \ $a, b \in C.$ Since we consider
the case of characteristic zero, the ideal of
identities with involution of $B$ is generated by ordinary
identities of $C,$ and it lies in $\Gamma.$ Indeed, take a multilinear
polynomial $f(y_1,\dots,y_k,z_{1},\dots,z_{n-k}) \in F\langle Y, Z
\rangle$ and consider a \mbox{$*$-homogeneous} substitution $y_1=(a_1,a_1),$
$\dots,$ $y_k=(a_k,a_k),$ $z_1=(a_{k+1},-a_{k+1}),$ $\dots,$
$z_{n-k}=(a_n,-a_n).$ We obtain
\begin{eqnarray*}
&&f((a_1,a_1),\dots,(a_k,a_k),(a_{k+1},-a_{k+1}),\dots,(a_n,-a_n))\\
&&=(f(a_1,\dots,a_k,a_{k+1},\dots,a_{n}), (-1)^{n-k}
\tilde{f}(a_1,\dots,a_k,a_{k+1},\dots,a_{n})),
\end{eqnarray*}
where $a_1,\dots,a_k,a_{k+1},\dots,a_{n} \in C$ are arbitrary.
Note that $f = 0$ is a polynomial
identity of an algebra with involution iff $\tilde{f} = 0$ is an
identity of the same algebra. Hence,
$f(y_1,\dots,y_k,z_{1},\dots,z_{n-k})=0$ is a $*$-identity of $B$ if
and only if $f(x_1,\dots,x_n),$ $\tilde{f}(x_1,\dots,x_n) \in
Id(C)=Id(A).$ Thus $Id^{*}(B) \subseteq Id^{*}(A).$
\hfill $\Box$

\section{The Kemer index of $*$-identity.}

Let us consider a finite dimensional $F$-algebra $A$ with involution.
It is well known (see, e.g., Theorems 3.4.3, 3.4.4 in \cite{GZbook}) that an analogue of the
Wedderburn-Malcev decomposition holds for algebras with involution.

\begin{lemma} \label{Pierce0} (\cite{GZbook})
Let $F$ be a field of zero characteristic. Then any finite dimensional
$F$-algebra with involution $A$ is isomorphic as $*$-algebra to an
$F$-algebra with involution of the form
\begin{equation} \label{matrix}
A'=C_1 \times \dots \times C_p \oplus J.
\end{equation}
Where the Jacobson radical $J=J(A')$ is a $*$-ideal, \
$B'=C_1 \times \dots \times C_p$ is a maximal semisimple
$*$-invariant subalgebra of $A',$ \ $C_l$ are $*$-simple algebras
($p  \ge 0$).
\end{lemma}

We can construct an algebra with involution which has the
Jacobson radical free in some sense. For a $*$-algebra $B$ (not necessarily without
unit) denote by $B^{\#}=B \oplus F \cdot 1_F$ the $*$-algebra
with adjoint unit $1_F$ \  ($(b + \alpha 1_F)^*=b^*+ \alpha 1_F$ \
for all $b \in B,$ \ $\alpha \in F$).

Given a finite dimensional $*$-algebra $A=B \oplus J(A)$ with a
maximal $*$-semisimple subalgebra $B$ and the Jacobson radical
$J(A)$ consider a $*$-invariant subalgebra $\widetilde{B}
\subseteq B.$ Take the free product $\widetilde{B}^{\#} *_F
F\langle X_q^{*} \rangle^{\#},$ where $F\langle X_q^{*} \rangle^{\#}$
is the free unitary \mbox{$*$-algebra}  of a rank
$q$. It is the $*$-algebra with involution naturally defined by
equalities $(u_1 \cdots u_s)^*=(u_s)^* \cdots (u_1)^*,$ where
$u_i \in \widetilde{B}^{\#} \bigcup F\langle X_q^{*}
\rangle^{\#}.$

Let us consider the subalgebra $\widetilde{B}(X_q^*)$ of
$\widetilde{B}^{\#} *_F F\langle X_q^{*} \rangle^{\#}$ generated
by $\widetilde{B} \cup F\langle X_q^{*} \rangle.$ It is
$*$-invariant. Denote by $(X_q^*)$ the two-sided $*$-ideal of
$\widetilde{B}(X_q^*)$ generated by the set of variables $X_q^*.$
Then $\widetilde{B}(X_q^*)=\widetilde{B} \oplus (X_q^*).$

Given a $*$T-ideal $\Gamma$ and a positive integer $s,$ we can consider the
quotient algebra
\begin{eqnarray} \label{FRad}
\mathcal{R}_{q,s}(\widetilde{B},\Gamma)=\widetilde{B}(X_q^*)/
(\Gamma(\widetilde{B}(X_q^*))+(X_q^*)^s).
\end{eqnarray}
Denote also
$\mathcal{R}_{q,s}(A)=\mathcal{R}_{q,s}(B,\mathrm{Id}^*(A)),$ where
$\Gamma=\mathrm{Id}^*(A),$ $\widetilde{B}=B.$

\begin{lemma} \label{Aqs}
For any $q, s \in \mathbb{N}$ and a $*$T-ideal $\Gamma \subseteq
\mathrm{Id}^*(A)$ the algebra
$\mathcal{R}_{q,s}(\widetilde{B},\Gamma)=\overline{B} \oplus
J(\mathcal{R}_{q,s}(\widetilde{B},\Gamma))$ is a finite dimensional
algebra with involution. Here $\overline{B} \cong \widetilde{B}$
is a maximal semisimple $*$-subalgebra of
$\mathcal{R}_{q,s}(\widetilde{B},\Gamma).$ The Jacobson radical
$J(\mathcal{R}_{q,s}(\widetilde{B},\Gamma))=(X_q^*)/
(\Gamma(\widetilde{B}(X_q^*))+(X_q^*)^s)$ is nilpotent of class
at most $s.$
$\mathrm{Id}^*(\mathcal{R}_{q,s}(\widetilde{B},\Gamma)) \supseteq
\Gamma.$

If $q \ge \mathrm{rk}(J(A))$, and $s$ is greater than
the class of nilpotency $\mathrm{nd}(A)$ of $J(A)$ then
$\mathrm{Id}^*(\mathcal{R}_{q,s}(A)) = \mathrm{Id}^*(A).$
\end{lemma}
\noindent {\bf Proof.} The ideal
$I=\Gamma(\widetilde{B}(X_q^*))+(X_q^*)^s$ is $*$-invariant. It is
clear that $I \subseteq (X_q^*).$  Then for the canonical
$*$-homomorphism $\psi:\widetilde{B}(X_q^*) \rightarrow
\mathcal{R}_{q,s}(\widetilde{B},\Gamma)$ we obtain
$\overline{B}=\psi(\widetilde{B}) \cong \widetilde{B},$ and
$\psi((X_q^*))=(X_q^*)/I$ is the maximal nilpotent ideal of the
algebra $\mathcal{R}_{q,s}(\widetilde{B},\Gamma)$ of class at
most $s.$ It is clear that $(X_q^*)/I$ is $*$-invariant and finite
dimensional. Hence,
$\mathcal{R}_{q,s}(\widetilde{B},\Gamma)=\overline{B} \oplus
\psi((X_q^*))$ is also a finite dimensional $*$-algebra with the
Jacobson radical
$J(\mathcal{R}_{q,s}(\widetilde{B},\Gamma))=\psi((X_q^*)).$ Also,
$\Gamma \subseteq
\mathrm{Id}^{*}(\mathcal{R}_{q,s}(\widetilde{B},\Gamma))$ for any
$q, s \in \mathbb{N}.$

Let us take $\Gamma=\mathrm{Id}^{*}(A),$ $\widetilde{B}=B.$ Consider
a generating set $\{r_1,\dots,r_\nu\}$ of the Jacobson radical
$J(A).$  If $q \ge \nu,$ and $s \ge
\mathrm{nd}(A)$ then the map $\varphi:\overline{x}_{i}=x_{i}+I
\mapsto r_{i}$ \ ($x_i \in X$ for $i=1,\dots,\nu$) can be extended
to a surjective $*$-homomorphism $\varphi:\mathcal{R}_{q,s}(A)
\rightarrow A$ assuming $\varphi(b+I)=b$ for all $b \in B.$
Therefore, $\Gamma = \mathrm{Id}^*(A) \supseteq
\mathrm{Id}^*(\mathcal{R}_{q,s}(A)).$ \hfill $\Box$

Observe that the construction of
$\mathcal{R}_{q,s}(\widetilde{B},\Gamma)$ can also be made
in terms of symmetric and skew-symmetric variables exchanging the
set of variables $X_q^*$ by the set $Y_q \cup Z_q=\{ y_1,\dots,
y_q, z_1,\dots, z_q \}.$

Over an algebraically closed field the Wedderburn-Malcev
decomposition (\ref{matrix}) of a finite dimensional $*$-algebra
can be described in more details.

\begin{lemma} \label{Pierce}
Let $F$ be an algebraically closed field. Any finite dimensional
$F$-algebra with involution $A$ is isomorphic as a $*$-algebra to an
$F$-algebra with involution $A'=C_1 \times \dots \times C_p \oplus
J,$ where $C_l$-s are $*$-simple algebras of the following types:
\begin{eqnarray} \label{list}
&& 1. \ (M_{k_l}(F),t) - \mbox{ the full matrix algebra with the transpose
involution } \ *=t, \nonumber \\ && 2. \ (M_{k_l}(F),s) - \mbox{
the full matrix algebra with the symplectic involution } *=s, \nonumber \\&& \ (k_l \in 2
\mathbb{Z}), \nonumber \\ && 3. \ (M_{k_l}(F) \times
M_{k_l}(F)^{op}, \bar{*}) - \ \mbox{ the direct product of the
  full matrix algebra } \nonumber \\ && \
  \mbox{ and its opposite algebra with the exchange
  involution } \  \bar{*} \mbox{ (see \cite{GZbook})}.
\end{eqnarray}
Moreover, $A'$ can be generated as a vector space by sets of its
symmetric elements $D_0,$ \ $U_0 \subseteq A'$, and skew-symmetric
elements $D_1,$ \ $U_1 \subseteq A'$ of the next form
\begin{eqnarray}
\label{basisD} &&D_0=\{ d^{(0)}_{l i_l j_l} = \varepsilon_l
d^{(0)}_{l i_l j_l}
 \varepsilon_l \ | \  \
(i_l, j_l) \in \mathcal{I}_l; \ 1 \le l \le p \}, \nonumber \\
&&D_1=\{ d^{(1)}_{l i_l j_l} = \varepsilon_l d^{(1)}_{l i_l j_l}
 \varepsilon_l \ | \  \ (i_l, j_l) \in \mathcal{J}_l; \ 1 \le l \le p \}, \\
 \label{basisU} &&U_0=\{
(\varepsilon_{l'} r \varepsilon_{l''} + \varepsilon_{l''} r^{*}
\varepsilon_{l'})/2| \ 1 \le l' \le l'' \le p+1; \ r \in J
\}\nonumber
\\ &&U_1=\{ (\varepsilon_{l'} r \varepsilon_{l''} -
\varepsilon_{l''} r^{*} \varepsilon_{l'})/2| \ 1 \le l' \le l''\le
p+1; \ r \in J \}.
\end{eqnarray}
Where \ \  in case  \quad $C_l=(M_{k_l}(F),t)$ \quad we set
\begin{eqnarray} \label{Clt}
&&d^{(0)}_{l i_l j_l} = E_{l i_l j_l} + E_{l j_l i_l}, \ \ 1 \le
i_l \le j_l \le k_l, \nonumber \\ &&d^{(1)}_{l i_l j_l} = E_{l i_l
j_l} - E_{l j_l i_l}, \ \ 1 \le i_l < j_l \le k_l;
\end{eqnarray}
in case  \quad $C_l=(M_{k_l}(F),s)$ \quad  we set
\begin{eqnarray} \label{Cls}
&&d^{(0)}_{l i_l j_l} = E_{l i_l j_l} + E_{l \frac{k_l}{2}+j_l
\frac{k_l}{2}+i_l}, \ \ 1 \le i_l, j_l \le \frac{k_l}{2},
\nonumber \\
&&d^{(1)}_{l i_l j_l} =E_{l i_l j_l} - E_{l
\frac{k_l}{2}+j_l \frac{k_l}{2}+i_l}, \ \  1 \le i_l, j_l \le
\frac{k_l}{2}, \nonumber \\
&&d^{(0)}_{l i_l \frac{k_l}{2}+j_l} =
E_{l i_l \frac{k_l}{2}+j_l} - E_{l j_l \frac{k_l}{2}+i_l}, \ \ 1
\le i_l < j_l \le \frac{k_l}{2}, \nonumber \\
&&d^{(1)}_{l i_l
\frac{k_l}{2}+j_l} = E_{l i_l \frac{k_l}{2}+j_l} + E_{l j_l
\frac{k_l}{2}+i_l} \ \  1 \le i_l \le j_l \le \frac{k_l}{2},
\nonumber \\
&&d^{(0)}_{l \frac{k_l}{2}+i_l j_l} = E_{l \frac{k_l}{2}+i_l j_l}
- E_{l \frac{k_l}{2}+j_l i_l},  \ \  1 \le i_l < j_l \le \frac{k_l}{2},
\nonumber
\\ &&d^{(1)}_{l \frac{k_l}{2}+i_l j_l} = E_{l \frac{k_l}{2}+i_l j_l}
+ E_{l \frac{k_l}{2}+j_l i_l}, \ \ 1 \le i_l \le j_l \le
\frac{k_l}{2};
\end{eqnarray}
in case  \quad $C_l=(M_{k_l}(F) \times M_{k_l}(F)^{op},\bar{*})$ \quad we set
\begin{eqnarray} \label{Clex}
&&d^{(0)}_{l i_l j_l} = (E_{l i_l j_l}, E_{l i_l j_l}), \ \ \ \  1
\le i_l, j_l \le k_l, \nonumber \\ &&d^{(1)}_{l i_l j_l} = (E_{l
i_l j_l}, - E_{l i_l j_l}), \ \ 1 \le i_l, j_l \le k_l.
\end{eqnarray}
Here $\varepsilon_l$ is the orthogonal central idempotent
of $B'=C_1 \times \dots \times C_p,$ corresponding to the unit
element of the $l$-th $*$-simple component $C_l$ of the algebra
$A'$ (for any $l=1,\dots,p$), $E_{l i_l j_l}$ are the matrix
units. Elements $r$ run on a set of elements of the Jacobson
radical $J=J(A)=\oplus_{l',l''=1}^{p+1} \varepsilon_{l'} J
\varepsilon_ {l''},$ where $\varepsilon_{p+1}$ is the adjoint
idempotent. All idempotents $\varepsilon_{l}$ ($l=1,\dots,p+1$)
are symmetric with respect to involution. The sets of pairs of
indices $\mathcal{I}_l,$ $\mathcal{J}_l$ are defined
by conditions (\ref{Clt}), (\ref{Cls}),
(\ref{Clex}).
\end{lemma}
\noindent {\bf Proof.} It is well known (see, e.g., Theorems
3.4.4 and 3.6.8 \cite{GZbook}) that any $*$-simple algebra is
isomorphic to one of the algebras of the list (\ref{list}). All
units $\varepsilon_l$ of $*$-simple components $C_l$ of $B'$ are
symmetric with respect to involution. Elements $d^{(0)}_{l
i_l j_l},$ \ $d^{(1)}_{l i_l j_l}$ defined by (\ref{Clt}),
(\ref{Cls}), (\ref{Clex}) form a \mbox{$*$-homogeneous} basis of the
corresponding $*$-simple component $C_l$ ($d^{(0)}_{l i_l
j_l}$ are symmetric, and \ $d^{(1)}_{l i_l j_l}$ are
skew-symmetric).

Any element $r \in J$ can be uniquely represented as a sum of
elements from the subspaces $\varepsilon_{l'} J \varepsilon_{l''}$
($l',l'' =1,\dots, p+1$) (the Peirce decomposition). All elements of
$\varepsilon_{l'} J \varepsilon_{l''}$ are sums of elements of
$U_0$ and $U_1$ for $1 \le l' \le l''\le p+1$, and differences of
elements of $U_0$ and $U_1$ if $l' > l''$ (for $r=v^*,$ \ $v \in
J$). Here elements of $U_0$ are symmetric, and the elements of
$U_1$ are skew-symmetric.  \hfill $\Box$

Observe that the condition for the base field $F$ to be algebraically closed
is necessary only for the description of the semi-simple part of $A$
((\ref{list}), (\ref{basisD}), (\ref{Clt}), (\ref{Cls}), (\ref{Clex})).
The statements concerning the Jacobson radical (for example, (\ref{basisU}))
is evidently true for any finite dimensional algebra with involution.
Moreover, considering \mbox{$*$-identities} we always can extend
statement of Lemma \ref{Pierce} to the case of arbitrary field
$F$ of zero characteristic.

\begin{definition}
An $F$-algebra $A$ is called representable if $A$ can be embedded
into some algebra $C$ that is finite dimensional over an extension
$\widetilde{F} \supseteq F$ of the base field $F.$
\end{definition}

\begin{lemma} \label{Repr}
Let $F$ be a field of zero characteristic. Any representable
$F$-algebra with involution $A$ is $*$PI-equivalent to some $F$-finite
dimensional algebra with involution $A'$ that satisfies
the claims of Lemma \ref{Pierce}.
\end{lemma}
\noindent {\bf Proof.}  We always can assume that the extension
$\widetilde{F} \supseteq F$  is algebraically
closed. Suppose that $A$ is isomorphic to an $F$-subalgebra
$\mathcal{B}$ of a finite dimensional $\widetilde{F}$-algebra
$\widetilde{\mathcal{B}}.$ The equality
$(\varphi(a))^*=\varphi(a^*)$ for any $a \in A$ ($\varphi:A
\rightarrow \mathcal{B}$ is the isomorphism of $F$-algebras) defines
involution on $\mathcal{B}.$
Let $\widetilde{U}$ be an $\widetilde{F}$-subalgebra of $*$-algebra
$\widetilde{\mathcal{B}} \times \widetilde{\mathcal{B}}^{op}$ with
the exchange involution $\bar{*}$ generated by $\{ (b,b^*) | b \in \mathcal{B} \}$
(here $*$ is the involution on $\mathcal{B}$ induced from $A$).
Then $\widetilde{U}=\{ \sum_{i=1}^{s} \alpha_i (b_i,b_i^*) |
\alpha_i \in \widetilde{F}, b_i \in \mathcal{B} \}$ is an $\widetilde{F}$-finite dimensional
$\bar{*}$-invariant subalgebra of
$\widetilde{\mathcal{B}} \times \widetilde{\mathcal{B}}^{op}.$
For algebras and identities over $F$ we have
$\mathrm{Id}^*(\widetilde{U})=\mathrm{Id}^{*}(\mathcal{B})=
\mathrm{Id}^{*}(A).$

Using Lemma \ref{Pierce}, over $\widetilde{F}$ we conclude that the
$*$-algebra $\widetilde{U}$ has the decomposition (\ref{matrix}),
where $\widetilde{C}_l$s are
isomorphic to  $(M_{k_l}(\widetilde{F}),t),$
$(M_{k_l}(\widetilde{F}),s),$ or $(M_{k_l}(\widetilde{F}) \times
M_{k_l}(\widetilde{F})^{op}, \bar{*}).$ It is clear from
(\ref{Clt}), (\ref{Cls}), and (\ref{Clex}) that
$\widetilde{C}_{l}=\widetilde{F} C_{l},$ where $C_l$ is
one of the algebras $(M_{k_l}(F),t),$
$(M_{k_l}(F),s),$ or $(M_{k_l}(F) \times M_{k_l}(F)^{op},
\bar{*})$ ($l=1,\dots,p$), respectively. Let us take $B=C_1 \times \dots \times
C_p,$ $\Gamma=\mathrm{Id}^*(A),$ $q=\dim_{\widetilde{F}}
J(\widetilde{U}),$ $s=\mathrm{nd}(\widetilde{U}).$ Then the
$F$-algebra $A'=\mathcal{R}_{q,s}(B,\Gamma)$ defined by
(\ref{FRad}) is a finite dimensional $F$-algebra with
involution. By Lemma \ref{Aqs} the algebra $A'$ satisfies the claims
of Lemma \ref{Pierce}, and $\mathrm{Id}^*(A') = \Gamma=
\mathrm{Id}^*(\widetilde{U}).$ \hfill $\Box$

In particular, any finite dimensional $F$-algebra with involution
$A$ can be naturally embedded to the $*$-algebra $A \otimes_F
\widetilde{F},$ which is finite dimensional over the algebraic
closure $\widetilde{F} \supseteq F.$ By this argument our study of
identities with involution can be reduced to finite
dimensional algebras specified in Lemma \ref{Pierce}.

We will assume without lost of generality that a finite
dimensional $*$-algebra $A$ over a field $F$ always satisfies
the claims of Lemma \ref{Pierce}. It means that a basis of $A$
can be always chosen in the set $D_0 \cup D_1 \cup U_0 \cup U_1.$
Notice that $D_0 \cup D_1$ is a \mbox{$*$-homogeneous} basis of the
semisimple part $B=C_1 \times \dots \times C_p$ of $A.$
Particularly, for multilinear $*$-polynomials it is enough to
consider only substitutions of $D_0 \cup U_0$ for
the symmetric variables, and $D_1 \cup U_1$ for
the skew-symmetric ones. A substitution of basic
elements from $D_0 \cup D_1 \cup U_0 \cup U_1$
(\ref{basisD}), (\ref{basisU}) for a multilinear $*$-polynomial
is called {\it elementary}. Elements of $D_0 \cup D_1$ are
called semisimple, whereas elements of  $U_0 \cup U_1$ are
called radical elements.

\begin{definition}
Let $A=B \oplus J$ be a finite dimensional algebra with involution
over a field $F,$ $B=B^{+} \oplus B^{-}$ a maximal semisimple
$*$-invariant subalgebra of $A,$ and $J(A)=J$ the Jacobson
radical of $A.$ Denote by $\mathrm{dims}_* A=(\dim B^{+},\dim
B^{-})$ dimensions of the symmetric and skew-symmetric parts of the
semisimple subalgebra $B,$ denote by $\mathrm{nd}(A)$ the
nilpotency class of the radical $J.$ Define the parameter of $A$
as \ $\mathrm{par}_*(A)=(\mathrm{dims}_* A;\mathrm{nd}(A)).$

The $4$-ple $\mathrm{cpar}_*(A)=(\mathrm{par}_*(A);\dim J(A))$
we call the complex parameter of $A.$
\end{definition}

Recall that $n$-tuples of numbers are ordered lexicographically.
It is clear that for any nonzero two-sided $*$-ideal $I \unlhd A$
of $A$ we have that $\mathrm{cpar}_*(A/I) < \mathrm{cpar}_*(A).$

Let us define principle numeric parameters of the ideal of
identities with involution of a finitely generated $*$-algebra.

Let $f=f(s_1,\dots,s_k,x_1,\dots,x_n) \in F\langle X^{*} \rangle$
be a polynomial linear in all variables $S=\{ s_1,
\dots, s_k\}$ ($S \subseteq Y$ or $S \subseteq Z$). We say
that $f$ is alternating in $S,$ if
$f(s_{\sigma(1)},\dots,s_{\sigma(k) },x_1,\dots,x_n)=(-1)^{\sigma}
f(s_1,\dots,s_k,x_1,\dots,x_n)$ holds for any permutation $\sigma
\in \mathrm{S}_k.$

For any polynomial with involution
$g(s_1,\dots,s_k,x_1,\dots,x_n)$ that is linear in $S=\{
s_1, \dots, s_k\}$ we construct a polynomial alternating in
$S$ by setting $$f(s_1,\dots,s_k,x_1,\dots,x_n)=
\mathcal{A}_{S}(g)=\sum_{\sigma \in \mathrm{S}_k} (-1)^{\sigma}
g(s_{\sigma(1)},\dots,s_{\sigma(k)},x_1,\dots,x_n).$$ The
corresponding mapping $\mathcal{A}_{S}$ is a linear transformation.
We call it the alternator. Any polynomial $f$ alternating in $S$
can be decomposed as $f=\sum_{i=1}^{m} \alpha_i
\mathcal{A}_{S}(u_i),$ where the $u_i$'s are monomials, \
$\alpha_i \in F.$ Properties of alternating polynomials with
involution are similar to that of ordinary polynomials (see,
e.g., \cite{Dren}, \cite{Kem1}, \cite{GZbook}). Note also
that a $*$-polynomial $f$ is alternating in $S$
iff its involution image $f^*$ is alternating in $S.$
Particularly, any alternator commutes with involution.

Given a pair $\overline{t}=(t^{+},t^{-}) \in \mathbb{N}_0^2$ we
say that a $*$-polynomial $f \in F\langle X^{*} \rangle$ has
$\overline{t}$-alternating variables
($f$ is $\overline{t}$-alternating) if $f(Y_0,Z_0,X)$ is linear in
$Y_0 \cup Z_0$ and $f$ is alternating in $Y_0 =
\{y_{i_1},\dots,y_{i_{t^{+}}}\} \subseteq Y,$ and $Z_0 =
\{z_{j_1},\dots,z_{j_{t^{-}}}\} \subseteq Z$ independently.

\begin{definition}
Fix $\overline{t}=(t^{+},t^{-}) \in \mathbb{N}_0^2.$
Suppose that $\tau_1, \dots, \tau_m \in \mathbb{N}_0^2$ are some
(possibly different) pairs satisfying the conditions
$\tau_j=(\tau_{j}^{+},\tau_{j}^{-})
> \overline{t}$ for $j=1, \dots, m.$
Let $f \in F\langle Y, Z \rangle$ be a multihomogeneous
$*$-polynomial in some symmetric and skew-symmetric
variables. Suppose that $f=f(S_1,\dots,S_{m+\mu};X)$ has $m$
collections of $\tau_j$-alternating variables
$S_j=Y_{j} \cup Z_{j}$ ($j=1,\dots,m$), \
and $\mu$ sets of $\overline{t}$-alternating variables $S_j=Y_{j} \cup Z_j$ \
($j=m+1,\dots,m+\mu$). We assume that all these sets are disjoint.
Then we say that $f$ is of the type
$(\overline{t};m;\mu)=(t^{+},t^{-};m;\mu).$ Here $Y_{j} \subseteq
Y$ are symmetric, and $Z_{j} \subseteq Z$ skew-symmetric
variables, and $|Y_{j}|=\tau_{j}^{+},$ $|Z_{j}|=\tau_{j}^{-}$ for
any $j=1,\dots,m,$ or $|Y_{j}|=t^{+},$ \ $|Z_{j}|=t^{-}$ \ for all
$j=m+1,\dots,m+\mu$ ($f$ is alternating in each $Y_j,$
$Z_j,$ \ $j+1,\dots,m+\mu$).
\end{definition}

Observe that a multihomogeneous polynomial $f$ of a type
$(\overline{t};m;\mu)$ is also of the type
$(\overline{t};m';\mu')$ for all $m' \le m,$ and $\mu' \le \mu.$
Particularly, any nontrivial multilinear $*$-polynomial in $Y \cup
Z$ of degree $m$ has the type $(0,0;m;\mu)$ for any $\mu \in
\mathbb{N}_0.$ Note also that multihomogeneous polynomials $f$ and
$f^*$ have always the same type.

\begin{definition} \label{defbeta}
Given a $*$T-ideal $\Gamma \unlhd F\langle Y, Z \rangle$ the
parameter $\beta(\Gamma)=(t^{+},t^{-})$ is the greatest
lexicographic pair $\overline{t}=(t^{+},t^{-}) \in \mathbb{N}_0^2$
such that for any $s \in \mathbb{N}$ there exists a $*$-polynomial
$f \notin \Gamma$ of the type $(\overline{t};0;s).$
\end{definition}
Any finitely generated PI-algebra $A$ satisfies the ordinary Capelli
identity of some order $d$ (\cite{Kem1}, \cite{GZbook}). Hence any polynomial
$f \in F\langle Y, Z \rangle$ of the type $(t_1,t_2;0;s)$ with
$t_1 \ge d$ or $t_2 \ge d$ ($s \in \mathbb{N}$)
belongs to $\mathrm{Id}^*(A).$ It means that the
parameter $\beta(\Gamma)$ is well defined for any proper
$*$T-ideal $\Gamma \unlhd F\langle Y, Z \rangle$ of identities
with involution of a finitely generated $*$-algebra. The next
parameter is also well defined.

\begin{definition} \label{defgamma}
Given a nonnegative integer $\mu$ let $\gamma(\Gamma;\mu)=s \in
\mathbb{N}$ be the smallest integer $s >0$ such that any
$*$-polynomial of type $(\beta(\Gamma);s;\mu)$ belongs to
$\Gamma.$

$\gamma(\Gamma;\mu)$ is a positive non-increasing function of
$\mu.$ Let us denote the limit of this function
$\gamma(\Gamma)=\lim \limits_{\mu \to \infty} \gamma(\Gamma;\mu)
\in \mathbb{N}.$

Then $\omega(\Gamma)$ is the smallest number $\widehat{\mu}$ such
that $\gamma(\Gamma;\mu)=\gamma(\Gamma)$ for any $\mu \ge
\widehat{\mu}$.
\end{definition}

\begin{definition}
We call by the Kemer index of a $*$T-ideal $\Gamma$ the
lexicographically ordered collection
$\mathrm{ind}_*(\Gamma)=(\beta(\Gamma); \gamma(\Gamma)).$
\end{definition}
Let us denote $\omega(A)=\omega(\mathrm{Id}^*(A)),$ \
$\gamma(A;\mu)=\gamma(\mathrm{Id}^*(A);\mu),$ \
$\mathrm{ind}_*(A)=\mathrm{ind}_*(\mathrm{Id}^*(A)),$ where $A$
is a finitely generated PI-algebra with involution.

It is clear that $A$ is nilpotent of class $s$ algebra with
involution if and only if
$\mathrm{ind}_*(A)=\mathrm{par}_*(A)=(0,0;s)$ ($\omega(A)=0$).

\begin{definition} \label{Kemerpolyn}
Given a nonnegative integer $\mu$ a nontrivial multihomogeneous
polynomial $f \in F\langle Y, Z \rangle$ is called $\mu$-boundary
polynomial for a $*$T-ideal $\Gamma$ if $f \notin \Gamma,$ and $f$
has the type $(\beta(\Gamma);\gamma(\Gamma)-1;\mu).$

Let us denote by $S_\mu(\Gamma)$ the set of all $\mu$-boundary
polynomials for $\Gamma$. Denote also \
$S_\mu(A)=S_\mu(\mathrm{Id}^*(A)),$ \
$K_\mu(\Gamma)=*T[S_\mu(\Gamma)],$ \ $K_{\mu,
A}=K_\mu(\mathrm{Id}^*(A))=*T[S_\mu(A)].$
\end{definition}

Any ideal of $*$-identities $\Gamma \ne F\langle Y, Z \rangle$  of
a finitely generated algebra has multilinear boundary polynomials
for all $\mu \in \mathbb{N}_0.$ Moreover, a polynomial $f$ belongs
to $S_{\mu}(\Gamma)$ along with its full multilinearization
$\widetilde{f},$ and its involution image $f^*.$

Note that the definitions of the Kemer index and boundary
polynomials for \mbox{$*$T-ideals} are adaptations of the corresponding
definitions for graded identities \cite{Svi}. Thus, the literal
repetition of the proofs of their properties for graded identities
\cite{Svi} also proves the same properties for identities with
involution.

\begin{lemma} \label{Prop1}
Given $*$T-ideals $\Gamma_1,$ $\Gamma_2 \subseteq F\langle Y, Z
\rangle$ admitting the Kemer indices, we have the next properties:
\begin{enumerate}
  \item If \  $\Gamma_1 \subseteq \Gamma_2$ \  then \  $\mathrm{ind}_*(\Gamma_1)
\ge \mathrm{ind}_*(\Gamma_2).$ \label{ind1}
  \item $\mathrm{ind}_*(\Gamma_1 \cap \Gamma_2)=\max
  \limits_{i=1,2} \mathrm{ind}_*(\Gamma_i).$ \label{dirsum}
  \item $\mathrm{ind}_*(A) \le \mathrm{par}_* (A)$ \  for any finite
dimensional  $*$-algebra $A.$ \label{bound}
\end{enumerate}
\end{lemma}

\begin{lemma} \label{subset}
For any $*$T-ideals $\Gamma_1,$ $\Gamma_2 \subseteq F\langle Y, Z
\rangle$ satisfying $\Gamma_1 \subseteq \Gamma_2$  one of the
following alternatives takes place:
\begin{enumerate}
  \item $\mathrm{ind}_*(\Gamma_1) = \mathrm{ind}_*(\Gamma_2),$ \
  and \ $S_\mu(\Gamma_1) \supseteq S_\mu(\Gamma_2)$ \ $\forall \mu
\in \mathbb{N}_0;$
  \item $\mathrm{ind}_*(\Gamma_1) > \mathrm{ind}_*(\Gamma_2),$ \
and \ $S_{\hat{\mu}}(\Gamma_1)
  \subseteq \Gamma_2$ for some $\hat{\mu} \in \mathbb{N}_0.$
\end{enumerate}
Moreover, in case $\Gamma_1 \subseteq \Gamma_2$ the conditions
\  $\mathrm{ind}_*(\Gamma_1) > \mathrm{ind}_*(\Gamma_2)$ and
$S_\mu(\Gamma_1) \subseteq \Gamma_2$ (correspondingly the
conditions $\mathrm{ind}_*(\Gamma_1) = \mathrm{ind}_*(\Gamma_2)$
and  $S_\mu(\Gamma_1) \supseteq S_\mu(\Gamma_2)$) are equivalent
for some $\mu \in \mathbb{N}_0.$
\end{lemma}

\begin{lemma} \label{Prop2}
Given $*$T-ideals $\Gamma,$ \ $\Gamma_1,$ $\dots,$ $\Gamma_{\rho}$
of $*$-identities of finitely generated PI-algebras with
involution we have:
\begin{enumerate}
  \item If \  $\mathrm{ind}_*(\Gamma_i)
< \mathrm{ind}_*(\Gamma)$ \ for all $i=1,\dots,\rho$ \  then there
exists $\widehat{\mu} \in \mathbb{N}_0$ such that $S_\mu(\Gamma)
\subseteq \bigcap \limits_{i=1} \limits^{\rho} \Gamma_i$ for any
$\mu \ge \widehat{\mu}.$ \label{ind2}

  \item If \  $\mathrm{ind}_* (\Gamma_i) =\kappa$ \  for any $i=1,\dots,\rho$
\  then for any $\mu \in \mathbb{N}_0$ holds
$$S_\mu(\bigcap \limits_{i=1} \limits^{\rho} \Gamma_i)=\bigcup
\limits_{i=1} \limits^{\rho} S_\mu(\Gamma_i), \qquad K_\mu(\bigcap
\limits_{i=1} \limits^{\rho}\Gamma_i)=\sum \limits_{i=1}
\limits^{\rho} K_\mu(\Gamma_i). $$ \label{dirsumKP}

  \item $\mathrm{ind}_* (\Gamma) > \mathrm{ind}_* (\Gamma+K_\mu(\Gamma))$
 \  for any $\mu \in \mathbb{N}_0.$  \label{addKP}
\end{enumerate}
\end{lemma}

\section{$*$PI-reduced algebras.}

Let us consider finite dimensional $*$-algebras which have the smallest
parameters for the same $*$-identities.

\begin{definition}
\qquad We say that a finite dimensional algebra $A$ with involution is \\
\mbox{$*$PI-reduced} if there do not exist finite dimensional $*$-algebras
$A_1,\dots,A_\varrho$ ($\varrho \in \mathbb{N}$) such that
$\bigcap\limits_{i=1}\limits^{\varrho} \mathrm{Id}^*(A_i) =
\mathrm{Id}^*(A),$ and $\mathrm{cpar}_*(A_i) < \mathrm{cpar}_*(A)$
for all $i=1,\dots,\varrho.$
\end{definition}

Similarly to the case of ordinary identities (\cite{BelRow},
\cite{Kem1}) and graded identities (\cite{AB}, \cite{Svi}) the
next natural facts are also take place for algebras with
involution.

\begin{lemma} \label{reduc}
Given a $*$PI-reduced algebra $A$ with the Wedderburn-Malcev
decomposition (\ref{matrix}) $A=(C_1 \times \cdots \times C_p)
\oplus J,$ we have $C_{\sigma(1)} J C_{\sigma(2)} J \cdots J
C_{\sigma(p)} \ne 0$ for some $\sigma \in \mathrm{S}_p.$ A
$*$PI-reduced algebra $A$ has no two proper two-sided $*$-ideals
$I_1, I_2 \lhd A$ satisfying $I_1 \cap I_2 = 0.$
\end{lemma}

Particularly, $\mathrm{nd}(A) \ge p$ always holds for a
$*$PI-reduced algebra $A$.

\begin{lemma} \label{decomp}
Any finite dimensional algebra with involution is $*$PI-equivalent
to a finite direct product of finite dimensional $*$PI-reduced
algebras.
\end{lemma}

Lemma \ref{decomp} along with Lemmas \ref{subset}, \ref{Prop2}
immediately implies the following.

\begin{lemma} \label{Smu}
A finite dimensional $*$-algebra $A$ is $*$PI-equivalent to
a direct product $\mathcal{O}(A) \times \mathcal{Y}(A)$ of finite
dimensional $*$-algebras $\mathcal{O}(A),$ $\mathcal{Y}(A)$
satisfying $\mathrm{ind}_*(A)=\mathrm{ind}_*(\mathcal{O}(A))
> \mathrm{ind}_*(\mathcal{Y}(A)).$ Moreover
$\mathcal{O}(A)=A_1 \times \dots \times A_\rho,$ where $A_i$ are
$*$PI-reduced, and
$\mathrm{ind}_*(A_i)=\mathrm{ind}_*(\mathcal{O}(A))$ for all
$i=1,\dots,\rho.$ There exists $\widehat{\mu} \in
\mathbb{N}_0$ such that
$S_\mu(A)=S_\mu(\mathcal{O}(A))=\bigcup_{i=1}^\rho S_\mu(A_i)
\subseteq \mathrm{Id}^*(\mathcal{Y}(A))$ holds for any $\mu \ge
\widehat{\mu}$.
\end{lemma}

\begin{definition} \label{sen}
$\mathcal{O}(A)$ is called the senior part of $A,$
$\mathcal{Y}(A)$ is called the junior part of $A.$ The algebras
$A_i$ ($i=1,\dots,\rho$) are called the senior components of $A.$
$\widehat{\mu}(A)$ is the minimal $\widehat{\mu} \in \mathbb{N}_0$
satisfying the conditions of Lemma \ref{Smu}.
\end{definition}

We have the next relations between the Kemer index and the parameters of
$*$PI-reduced algebras.

\begin{lemma} \label{ind-simple}
Given a $*$PI-reduced algebra $A$ we have
$\beta(A)=\mathrm{dims}_* A.$ Any $*$-simple finite dimensional
algebra $C$ is $*$PI-reduced, and
$\mathrm{ind}_*(C)=\mathrm{par}_*(C)=(t_1,t_2;1).$
\end{lemma}
\noindent {\bf Proof.} For a nilpotent algebra $A$ the assertion
of lemma is trivial. Assume that $A$ is a non-nilpotent
$*$PI-reduced algebra with $\mathrm{dims}_* A=(t_1,t_2).$ From
Lemma \ref{Prop1} it follows that $\beta(A) \le \mathrm{dims}_*
A.$ Thus it is enough to construct for any $\hat{s} \in
\mathbb{N}$ a $*$-polynomial of the type $(t_1,t_2;0;\hat{s})$
which is not an identity with involution of $A.$

We can assume that $A$ has decomposition (\ref{matrix})
specified by Lemmas \ref{Pierce0}, \ref{Pierce}. Consider any
$*$-simple component $C_l$ ($l=1,\dots,p$), and take $\hat{s}$ sets of
distinct symmetric variables $Y_{l, m}=\{ y_{l, m, (i_l j_l)} \in
Y | (i_l,j_l) \in \mathcal{I}_l \}$ corresponding to the symmetric
basic elements $d^{(0)}_{l i_l j_l},$ and $\hat{s}$ sets of
distinct skew-symmetric variables $Z_{l, m}=\{z_{l, m, (i_l j_l)}
\in Z | (i_l,j_l) \in \mathcal{J}_l \}$ corresponding to the
skew-symmetric basic elements $d^{(1)}_{l i_l j_l}$ (see
(\ref{basisD})), where $m=1,\dots,\hat{s}$. Consider also the set
of ordinary variables (neither symmetric nor skew-symmetric)
$X_{l}=\{ x_{l, (i_l j_l,i'_l j'_l)} |
1 \le i_l, j_l, i'_l, j'_l \le k_l \} \subseteq X.$
We say that the variable $x_{l, (j_1 i_2,j_2 i_1)}$ connects the
variables $y_{l, m, (i_1 j_1)}$ and $y_{l, m, (i_2 j_2)}$ (or the
variables $z_{l, m, (i_1 j_1)}$ and $z_{l, m, (i_2 j_2)}$).

Let $l,$ $m$ be fixed, we consider the $*$-monomial $w_{l, m
}(Y_{l, m},Z_{l, m},X_{l})$ that is the product of all variables
$y_{l, m, (i_l j_l)},$ and all variables $z_{l, m,(i_l j_l)}$
connected by the variables $x_{l, (i j,i' j')},$ correspondingly ($l=1,
\dots, p,$ \ $m=1,\dots,\hat{s}$). For example, if
$C_l=(M_{3}(F),t)$ then we obtain
\begin{eqnarray*}
&&w_{l, m}= y_{l, m,(11)} x_{l, (11,21)} y_{l, m, (12)} x_{l, (21,31)}
y_{l, m, (13)} x_{l, (32,21)} y_{l, m, (22)} x_{l, (22,32)}  y_{l, m,
(23)} \\ && x_{l, (33,32)} y_{l, m, (33)} x_{l, (31,23)} z_{l, m, (12)}
x_{l, (21,31)} z_{l, m, (13)} x_{l, (32,31)} z_{l, m, (23)}.
\end{eqnarray*}
The unit element of $C_l$ can be decomposed as
$\varepsilon_l=\sum_{s_l=1}^{k_l} E_{l s_l s_l}$ in cases
(\ref{Clt}), (\ref{Cls}), and $\varepsilon_l=\sum_{s_l=1}^{k_l}
(E_{l s_l s_l},E_{l s_l s_l})$ in case (\ref{Clex}).
Let us denote also $e^{(l)}_{(i_l j_l,i'_l j'_l)}=E_{l i_l j_l}$
for all $i'_l, j'_l$ for algebras $C_l$ of the types (\ref{Clt}), (\ref{Cls}), or
$e^{(l)}_{(i_l j_l,i'_l j'_l)}=(E_{l i_l j_l},E_{l i'_l j'_l})$
for the algebra $C_l$ of the form (\ref{Clex}).
Thus by
Lemma \ref{reduc}, we can assume that $A$ contains an element
\begin{equation} \label{ErE}
e^{(1)}_{s_1 s_1} r_1 e^{(2)}_{s_2 s_2} \cdots e^{(p-1)}_{s_{p-1}
s_{p-1} } r_{p-1} e^{(p)}_{s_{p} s_{p} } \ne 0,
\end{equation}
where $r_l \in J$ are some radical elements, and $e^{(l)}_{i_l
j_l}=e^{(l)}_{(i_l j_l,i_l j_l)}.$

Consider the $*$-monomial $W_l=x_{l, (s_l t_l,t'_l s_l)}\cdot \bigl(
\prod_{m=1}^{\hat{s}}  (x_{l, (t_l 1,1 t'_l)} \cdot w_{l, m}) \bigr)
 \cdot x'_{l, (t_l s_l,s_l t'_l)},$ where $s_l$ are given by
(\ref{ErE}), and $t_l,$ $t'_l$ are chosen to connect the word $w_{l, m}$
with $w_{l, m+1}.$ Let us denote $Y_{(m)}= \bigcup_{l=1}^{p}
Y_{l,m},$ and $Z_{(m)}= \bigcup_{l=1}^{p} Z_{l,m}.$ It is clear
that $|Y_{(m)}|=t_1,$ $|Z_{(m)}|=t_2$ for any $m=1,\dots,\hat{s}.$
Then the polynomial
\begin{eqnarray} \label{polyn2}
f=\bigl( \prod_{m=1}^{\hat{s}} \mathcal{A}_{Y_{(m)}}
\mathcal{A}_{Z_{(m)}} \bigr) (W_1 \tilde{x}_1 W_2 \tilde{x}_2
\cdots \tilde{x}_{p-1} W_p)
\end{eqnarray}
is alternating in $Y_{m} \subseteq Y$ and
$Z_{m} \subseteq Z$ \ for all $m=1,\dots,\hat{s}.$ Here the
variables $X_{l} \subseteq X,$ \
$X'_{l}=\{ x'_{l, (i_l j_l,i'_l j'_l)} | 1 \le i_l, j_l,
i'_l, j'_l \le k_l \}  \subseteq X$ and $\tilde{x}_l
\in X$ are not \mbox{$*$-homogeneous} ($l=1,\dots,p$). Consider the substitution
for $f$ as follows
\begin{eqnarray} \label{s-sim}
&&y_{l, m, (i_l j_l)} = d^{(0)}_{l i_l j_l}, \  ((i_l,j_l) \in
\mathcal{I}_l); \qquad z_{l, m, (i_l j_l)} = d^{(1)}_{l i_l j_l},
\ ((i_l,j_l) \in \mathcal{J}_l); \nonumber \\ &&x_{l, (i_l
j_l,i'_l j'_l)}=e^{(l)}_{(i_l j_l,i'_l j'_l)}, \
( 1 \le i_l,j_l \le k_l); \qquad
\tilde{x}_{q}=\varepsilon_q r_q \varepsilon_{q+1}; \nonumber
\\ &&x'_{l, (i_l
j_l,i'_l j'_l)}=(E_{l i_l j_l},(-1)^{k_l} E_{l i'_l j'_l}) \
\mbox{ if } C_l=(M_{k_l}(F) \times M_{k_l}(F)^{op},\bar{*}) \
(\mbox{ see } (\ref{Clex})), \nonumber
\\&&x'_{l, (i_l
j_l,i'_l j'_l)}=e^{(l)}_{i_l j_l} \ \mbox{ otherwise, }
\  ( 1 \le i_l,j_l,i'_l,j'_l \le k_l); \nonumber
\\&&l=1,\dots,p; \ \ q=1,\dots, p-1; \ \ m=1,\dots,\hat{s}.
\end{eqnarray}
The elements from the simple component $C_l$ cannot be
substituted for the variables of $W_q,$ \  $l \ne q,$ with a nonzero result. Therefore,
the result of the substitution (\ref{s-sim}) applied to $f$ coincides with
the result of the same substitution applied to the polynomial $$f'=\Bigl(
( \prod_{m=1}^{\hat{s}} \mathcal{A}_{Y_{1, m}} \mathcal{A}_{Z_{1,
m}}) W_1  \Bigr) \tilde{x}_1  \cdots \tilde{x}_{p-1} \Bigl(
(\prod_{m=1}^{\hat{s}} \mathcal{A}_{Y_{p, m}} \mathcal{A}_{Z_{p,
m}}) W_p \Bigr).$$ For any $m$ the set $Y_{l, m}$ contains at most
one variable $y_{l, m, (i_l j_l)}$ for the same pair $(i_l,j_l) \in
\mathcal{I}_l$ (correspondingly, the set $Z_{m}$ contains at most one
variable $z_{l, m, (i_l j_l)}$ for the same pair $(i_l,j_l) \in
\mathcal{J}_l$). Since we fix positions for $(i_l, j_l)$ in the
direct (and in the opposite if necessary) component of $C_l$
uniquely by $x_{l, (i'_l j'_l,i''_l j''_l)}=
e^{(l)}_{(i'_l j'_l,i''_l j''_l)}$
then the substitution
(\ref{s-sim}) applied to the polynomial $f^{(l)}=(\prod_{m=1}^{\hat{s}}
\mathcal{A}_{Y_{l, m}} \mathcal{A}_{Z_{l, m}}) W_l$ \ \  gives the
result $2^{\xi_l} e^{(l)}_{s_l s_l},$ \ $\xi_l \in \mathbb{N}_0.$
Hence the result of the substitution
(\ref{s-sim}) to the polynomial $f$ gives the nonzero element
(\ref{ErE}) of $A.$ Therefore, we obtain $f \notin
\mathrm{Id}^*(A).$

Consider the exchange (\ref{free}) of variables $\bigcup_{l=1}^{p}
\bigl(X_{l} \cup X'_{l} \cup \{\tilde{x}_l \}\bigr)$ by the sums of
corresponding symmetric and skew-symmetric variables. Then at
least one multihomogeneous component $\tilde{f} \in F\langle Y, Z
\rangle$ of $f$ is also not a $*$-identity of $A.$ Recall that $\tilde{f}$ is
alternating in $Y_{(m)},$ and $Z_{(m)}$
($m=1,\dots,\hat{s}$). Thus $\tilde{f}$ is the required
polynomial.

Notice that the condition of $*$PI-reducibility of $A$ is used in
the proof only to find a non-zero element (\ref{ErE}) in $A.$
Thus in the case of a $*$-simple algebra $A$ ($p=1$) this
condition is not necessary. In this case
we can take $e^{(1)}_{1 1}$ instead of (\ref{ErE}), and we also obtain
$\beta(A)=\mathrm{dims}_* A.$ A finite dimensional $*$-simple
algebra is semisimple, $\mathrm{nd}(A)=1,$ and
$\mathrm{cpar}_*(A)=(\mathrm{dims}_* A;1;0).$ Taking into account
that $\mathrm{ind}_*(A) \le \mathrm{par}_*(A)=(\beta(A);1)$ (Lemma
\ref{Prop1}), and $\gamma(A)>0$ we obtain
$\mathrm{ind}_*(A) = \mathrm{par}_*(A)$ for a $*$-simple algebra.
By Lemma \ref{Prop1} the
conditions $\dim J(A)=0,$ and $\mathrm{ind}_*(A) =
\mathrm{par}_*(A)$ imply also that $A$ is $*$PI-reduced. \hfill
$\Box$

\begin{definition}
Assume that a finite dimensional algebra with involution $A$
satisfies the claims of Lemma \ref{Pierce}. An elementary
substitution $(a_1,\dots,a_n)$ of $*$-homogeneous elements of $A$
(namely, $a_i \in D_0 \cup D_1 \cup U_0 \cup U_1 \subseteq A$
(\ref{basisD}), (\ref{basisU})) is called incomplete if there
exists $j \in \{ 1,\dots,p \}$ such that $$\{a_1,\dots,a_n\} \cap \bigl(C_j
\oplus_{l=1}^{p+1} (\varepsilon_j J \varepsilon_l + \varepsilon_l
J \varepsilon_j) \bigr) = \emptyset.$$ Otherwise, the substitution
$(a_1,\dots,a_n)$ is called complete.
\end{definition}

\begin{definition}
An elementary substitution $(a_1,\dots,a_n) \in A^n$ is
called thin if it contains less than $\mathrm{nd}(A)-1$
radical elements (not necessarily distinct).
\end{definition}

\begin{definition}
We say that a multilinear $*$-polynomial
$f(y_1,\dots,y_{n_1},z_1,\dots,z_{n_2}) \in F\langle Y, Z \rangle$
is exact for a finite dimensional $*$-algebra $A$ if
$f(a_1,\dots,a_n)=0$ holds in $A$ for any thin or incomplete
substitution $(a_1,\dots,a_n) \in A^n$ ($n=n_1+n_2$).
\end{definition}

\begin{lemma} \label{Exact2}
If $A$ is a $*$PI-reduced algebra then any multilinear polynomial
of the type $(\mathrm{dims}_* A;\mathrm{nd}(A)-1;0)$ is exact for
$A.$
\end{lemma}
\noindent {\bf Proof.} It is clear that a $*$PI-reduced algebra is
either not semisimple or $*$-simple. For a $*$-simple finite
dimensional algebra any multilinear $*$-polynomial can be assumed
exact. A multilinear polynomial of type
$(0,0;\mathrm{nd}(A)-1;0)$ has the full degree greater or equal to
$(\mathrm{nd}(A)-1).$ Hence it is assumed to be exact for a
nilpotent algebra $A.$

Suppose that $A$ is neither nilpotent nor semisimple, and $f \in
F\langle Y, Z \rangle$ is a multilinear $*$-polynomial of the type
$(\mathrm{dims}_* A;\mathrm{nd}(A)-1;0).$ Under a thin
substitution at least one of $\tilde{s}=\mathrm{nd}(A)-1$
collections of $\tau_{j}$-alternating variables of $f$ will be
completely replaced by semisimple elements. Since
$\tau_{j}>\mathrm{dims}_* A$ for any $j=1,\dots,\tilde{s}$ then
the result of the substitution will be zero.

By Lemma \ref{Pierce}, it is clear that $\dim_F C_{l}^{+}
> 0$ for any $l=1,\dots,p$. Therefore, an incomplete
substitution cannot contain all semisimple symmetric elements
from $D_0$ (\ref{basisD}). Taking into account the conditions
$\tau_{j}>\mathrm{dims}_* A,$ \  $j=1,\dots,\tilde{s},$ we
obtain that at least two variables of every collection of
$\tau_j$-alternating variables of $f$ must be replaced by
radical elements, otherwise the result of the substitution will be
zero. Thus, the result of an incomplete substitution to
the polynomial $f$ is zero. \hfill $\Box$

\begin{lemma} \label{Exact1}
Any nonzero $*$PI-reduced algebra $A$ has an exact polynomial,
that is not a $*$-identity of $A.$
\end{lemma}
\noindent {\bf Proof.} For a nilpotent algebra $A$ the assertion
follows from Lemma \ref{Exact2}. Suppose that $A$ is a
non-nilpotent algebra satisfying the claims of Lemma
\ref{Pierce}. Consider its subalgebras $A_i=(\prod
\limits_{\mathop{1 \le j \le p}\limits_{\scriptstyle j \ne i}} C_j
) \oplus J(A)$ for all $i=1,\dots, p$. Take $q=\dim_F J(A),$ \
$s=\mathrm{nd}(A)-1.$ Then by Lemma \ref{Aqs} \ $\widetilde{A}=A_1
\times \dots \times A_p \times \mathcal{R}_{q,s}(A)$ is a finite
dimensional $*$-algebra satisfying $\mathrm{Id}^*(A) \subseteq
\mathrm{Id}^*(\widetilde{A}).$  Let $\{r_1,\dots,r_q\} \subseteq
U_0 \cup U_1$ be a $*$-homogeneous basis of $J(A)$ of the form
(\ref{basisU}). Consider the map $\varphi$ defined by
$\varphi(y_{j})=r_j$ if $r_j$ is symmetric, \ $\varphi(z_{j})=r_j$
if $r_j$ is skew-symmetric (for all $i=1,\dots,q$), and
$\varphi(b)=b$ for any $b \in B.$ Then $\varphi$ can be extended to a surjective
$*$-homomorphism $\varphi: B(Y_q,Z_q) \rightarrow A.$ It follows
that any multilinear polynomial $f \in
\mathrm{Id}^*(\mathcal{R}_{q,s}(A))$ is turned into zero under any thin
substitution. It is also clear that any incomplete substitution in
a multilinear $*$-polynomial $f \in \mathrm{Id}^*(\times_{i=1}^{p}
A_i)$ yields zero.  Therefore, any multilinear polynomial
$f \in \mathrm{Id}^*(\widetilde{A})$ is exact for $A.$ Remark that
$\mathrm{cpar}_*(A_i) < \mathrm{cpar}_*(A)$ ($1 \le i \le p$), and
$\mathrm{cpar}_*(\mathcal{R}_{q,s}(A)) < \mathrm{cpar}_*(A).$
Since $A$ is a $*$PI-reduced algebra then $\mathrm{Id}^*(A)
\subsetneqq \mathrm{Id}^*(\widetilde{A}).$ Any multilinear
polynomial $f$ such that $f \in \mathrm{Id}^*(\widetilde{A}),$ and
$f \notin \mathrm{Id}^*(A)$ satisfies the assertion of the lemma.
\hfill $\Box$

\begin{lemma} \label{Gammasub}
Let $A$ be a finite dimensional $*$-algebra, $h$ an
exact $*$-polynomial for $A,$  and \  $\bar{a} \in A^n$ a complete
substitution in $h$ containing exactly $\tilde{s}=\mathrm{nd}(A)-1$
radical elements. Then for any $\mu
\in \mathbb{N}_0$ there exist a $*$-polynomial
$h_\mu \in *T[h],$ and
an elementary substitution $\bar{u}$ from the algebra $A$ into $h_\mu$ such that:
\begin{enumerate}
 \item $h_\mu(\mathcal{Z}_1,\dots,\mathcal{Z}_{\tilde{s}+\mu},
\mathcal{X})$ is
$\tau_j$-alternating in any set $\mathcal{Z}_j$ with $\tau_j >
\beta=\mathrm{dims}_* A$ for all $j=1,\dots,\tilde{s}$, and
is $\beta$-alternating in any $\mathcal{Z}_j$ for
$j=\tilde{s}+1,\dots,\tilde{s}+\mu$ (all the sets
$\mathcal{Z}_j,$ $\mathcal{X} \subseteq (Y \cup Z)$ are disjoint),
 \item $h_\mu(\bar{u})= \alpha h(\bar{a})$ for some $\alpha \in F,$ $\alpha \ne 0,$
 \item all
variables from $\mathcal{X}$ are
replaced by semisimple elements.
\end{enumerate}
\end{lemma}
\noindent {\bf Proof.}
If $h(\bar{a})=0$ then the assertion of lemma is trivial. It is sufficient
to take any consequence $h_\mu(\mathcal{Z}_1,\dots,\mathcal{Z}_{\tilde{s}+\mu}) \in *T[h]$
that is alternating in $\mathcal{Z}_j$ as required (we assume here
that $\mathcal{X}=\emptyset$), and replace the variables of the alternating set $Z_j$
by equal elements. Particularly, from conditions $\mathrm{nd}(A)=1,$ \ $p \ge 2$ it follows that $h(\bar{a})=0.$

Assume that $h(\bar{a}) \ne 0.$
Let us consider the case $\mathrm{nd}(A)>1,$ \ $p
\ge 2$ in decomposition (\ref{matrix}) of the algebra $A.$ We can suppose
for simplicity that the substitution $\bar{a}$ has the form
\quad $\bar{a}=(r_{1, (l'_1,l''_1)},\dots,r_{\tilde{s}, (l'_{\tilde{s}},l''_{\tilde{s}})},
b_1,\dots,b_{n-\tilde{s}}),$ \ where
$1 \le l'_s,l''_s \le p+1,$ \ $\{1,\dots,p\} \subseteq \{l'_s, l''_s
\ | 1 \le s \le q \}$ for some $q \le \tilde{s},$ \ $l'_s \ne l''_s$ for any $s=1,\dots,q.$ \
Here we have elements $b_1,\dots,$ $b_{n-\tilde{s}} \in D_0 \cup D_1$ (\ref{basisD}),
and $r_{s, (l'_s,l''_s)} \in U_0 \cup U_1$ (\ref{basisU}), where $r_{s, (l'_s,l''_s)} =
(\varepsilon_{l'_s} r_s \varepsilon_{l''_s} + \varepsilon_{l''_s} r_s^{*}
\varepsilon_{l'_s})/2$ \  if
$r_{s, (l'_s,l''_s)}$ is a symmetric element of $A$, or $r_{s, (l'_s,l''_s)}=
(\varepsilon_{l'_s} r_s \varepsilon_{l''_s} - \varepsilon_{l''_s} r_s^{*}
\varepsilon_{l'_s})/2$ \
if $r_{s, (l'_s,l''_s)}$ is skew-symmetric ($s=1,\dots,\tilde{s}$).
Particularly, without lost of generality we can consider the substitution
$\bar{a}=(r_{1, ({\bf 1},{\bf 2})},r_{2, ({\bf 3},l_2'')},\dots,r_{q, (l'_{q},{\bf p})},\dots,
b_1,\dots,b_{n-\tilde{s}}).$

Let us take for any $l=1,\dots,p$ the monomials
$w_{l, m}(Y_{l, m},Z_{l, m},X_{l})$ specified in Lemma \ref{ind-simple}, where
$m=1,\dots,\tilde{s}+\mu.$ Construct the $*$-polynomial
$\tilde{f}_l(Y_{l},Z_{l},\widetilde{X}_{l})=\sum_{s_l=1}^{k_l}
\bigl( x_{l, (s_l t_l,t'_l s_l)}\cdot \bigl(
\prod_{m=1}^{\tilde{s}+\mu}  (x_{l, (t_l 1,1 t'_l)} \cdot w_{l, m}) \bigr)
 \cdot x'_{l, (t_l s_l,s_l t'_l)} \bigr),$ where $t_l,$ $t'_l$ connect the
words $w_{l, m}$ and $w_{l, m+1}$ like in Lemma \ref{ind-simple},
and $Y_{l}=\bigcup_{m=1}^{\tilde{s}+\mu} Y_{l, m},$ \
$Z_{l}=\bigcup_{m=1}^{\tilde{s}+\mu} Z_{l, m}.$ Make the replacement
$x_{l, (i_l j_l, i'_l j'_l)}=(\delta_1 \
\tilde{y}_{l, (i_l j_l)}+ \delta_2 \ \tilde{z}_{l, (i_l j_l)} +
\delta_3 \ \tilde{y}_{l, (i'_l j'_l)}
+\delta_4 \ \tilde{z}_{l, (i'_l j'_l)})/2,$ \
$x'_{l, (i_l j_l, i'_l j'_l)}=(\delta_5 \ \tilde{y}'_{l, (i_l j_l)}+
\delta_6 \ \tilde{z}'_{l, (i_l j_l)}+ \delta_7 \
\tilde{y}'_{l, (i'_l j'_l)}+ \delta_8 \
\tilde{z}'_{l, (i'_l j'_l)})/2$
of the variables $x_{l, (i_l j_l,i'_l j'_l)},
x'_{l, (i_l j_l,i'_l j'_l)}  \in \widetilde{X}_{l}$
by the combination of corresponding symmetric and skew-symmetric variables
$\widetilde{Y}_{l}=\{ \tilde{y}_{l, (i_l j_l)},
\tilde{y}'_{l, (i_l j_l)} | 1 \le i_l, j_l \le k_l  \},$ \
$\widetilde{Z}_{l}=\{ \tilde{z}_{l, (i_l j_l)},
\tilde{z}'_{l, (i_l j_l)} | 1 \le i_l, j_l \le k_l  \}$
pairwise different and disjoint with $Y_{l} \cup Z_{l},$ \
$\delta_c \in \{ 0,1,-1\}$ ($c=1,\dots,8$).
Let us denote by $f_l(Y_{l},Z_{l},\widetilde{Y}_{l},\widetilde{Z}_{l})=
(\tilde{f}_l+\tilde{f}_l^{*})/2 \in F \langle Y, Z \rangle$ the symmetric part
of the $*$-polynomial $\tilde{f}_l$ after this replacement.
Consider the polynomials
\begin{eqnarray*}
h'(\mathcal{Z}_1,\dots,\mathcal{Z}_{\tilde{s}+\mu},
\mathcal{X})=
h \Bigl(f_1 \circ f_2 \circ \zeta_1, f_3 \circ \zeta_2, \dots, f_{p} \circ \zeta_q,
\zeta_{q+1},\dots,\zeta_{\tilde{s}}, \
\tilde{x}_1,\dots,\tilde{x}_{n-\tilde{s}} \Bigr); \\
\mbox{and } \qquad h_\mu(\mathcal{Z}_1,\dots,\mathcal{Z}_{\tilde{s}+\mu},
\mathcal{X})= \Bigl( \prod_{m=1}^{\tilde{s}+\mu}
(\mathcal{A}_{\mathcal{Z}^{(0)}_m}  \cdot \mathcal{A}_{\mathcal{Z}^{(1)}_m}) \Bigr)
h'(\mathcal{Z}_1,\dots,\mathcal{Z}_{\tilde{s}+\mu},
\mathcal{X}).
\end{eqnarray*}
Here $\mathcal{Z}_m=\bigcup_{l=1}^p (Y_{l, m} \cup Z_{l, m}) \cup \{ \zeta_m \}$
if $m=1,\dots,\tilde{s},$ and $\mathcal{Z}_m=\bigcup_{l=1}^p (Y_{l, m} \cup Z_{l, m})$
for $m=\tilde{s}+1,\dots,\tilde{s}+\mu,$ \
$\mathcal{X}=\bigcup_{l=1}^p (\widetilde{Y}_{l} \cup \widetilde{Z}_{l}) \cup
\{ \tilde{x}_1,\dots,\tilde{x}_{n-\tilde{s}} \}.$
Here the variable $\zeta_s$ is symmetric if the element $r_{s, (l'_s,l''_s)}$
of the substitution $\bar{a}$ is symmetric, and $\zeta_s$ is skew-symmetric
when $r_{s, (l'_s,l''_s)}$ is skew-symmetric. Similarly, a variable $\tilde{x}_i$
is symmetric or skew-symmetric according to the respective element $b_i$ of
the substitution $\bar{a}.$ Notice that $\mathcal{Z}^{(0)}_m$
is a subset of all symmetric variables of $\mathcal{Z}_m,$
and $\mathcal{Z}^{(1)}_m$ is a subset of all skew-symmetric
variables of $\mathcal{Z}_m.$
Particularly, $\mathcal{Z}^{(0)}_m=\bigcup_{l=1}^p Y_{l, m},$ \
$\mathcal{Z}^{(1)}_m=\bigcup_{l=1}^p Z_{l, m}$ for any
$m=\tilde{s}+1,\dots,\tilde{s}+\mu.$
Since any polynomial $f_l$ is a symmetric element of the free algebra
with involution then $f_l \circ \zeta_s$ is also symmetric (skew-symmetric) as soon as
$\zeta_s$ is symmetric (skew-symmetric). Thus, the polynomial
$h_{\mu} \in *T[h]$ is well defined.

It is clear that the polynomial $h_{\mu}$ is
$\tau_j$-alternating in $\mathcal{Z}_j$ with $\tau_j >
\beta=\mathrm{dims}_* A$ for all $j=1,\dots,\tilde{s}$, and
is $\beta$-alternating in $\mathcal{Z}_j$ for
$j=\tilde{s}+1,\dots,\tilde{s}+\mu.$

Consider the next substitution of the polynomial $h_{\mu}$
\begin{eqnarray} \label{substK22}
&&y_{l, m, (i_l j_l)} = d^{(0)}_{l i_l j_l}, \  ((i_l,j_l) \in
\mathcal{I}_l); \qquad z_{l, m, (i_l j_l)} = d^{(1)}_{l i_l j_l},
\ ((i_l,j_l) \in \mathcal{J}_l); \nonumber \\
&&\zeta_{s}=a_s=r_{s, (l'_s, l''_s)}; \qquad \qquad \tilde{x}_{n'}=
a_{n'+\tilde{s}}=b_{n'};
\\ &&l=1,\dots,p; \ \ m=1,\dots,\tilde{s}+\mu;
\ \  1 \le s \le \tilde{s}; \ \  1
\le n' \le n-\tilde{s}. \nonumber
\end{eqnarray}
The elementary substitution of the variables of the sets $\widetilde{Y}_{l},$
$\widetilde{Z}_{l},$ and the coefficients $\delta_c \in \{0, 1, -1 \}$
are chosen to guarantee $x_{l, (i_l j_l,i'_l j'_l)}=
e^{(l)}_{(i_l j_l,i'_l j'_l)}$ ($1 \le i_l, j_l, i'_l, j'_l \le k_l$,
$1 \le l \le p$); and \  $x'_{l, (i_l j_l,i'_l j'_l)}=
e^{(l)}_{i_l j_l}$ ($1 \le i_l, j_l, i'_l, j'_l \le k_l$)
if $C_l=(M_{k_l}(F),*)$ with $* \in \{ t, s \}$ (\ref{Clt}), (\ref{Cls}),
or $x'_{l, (i_l j_l,i'_l j'_l)}=(E_{l i_l j_l},(-1)^{k_l} E_{l i'_l j'_l})$
if $C_l=(M_{k_l}(F) \times M_{k_l}(F)^{op},\bar{*})$ (\ref{Clex}).

Due to the substitution of the variables of
$\widetilde{Y}_{l} \cup \widetilde{Z}_{l},$
the polynomial $f_l$ can contain only
elements of the simple component $C_l$ or elements
$r_{s, (l'_s,l''_s)}$ with $l'_s=l''_s=l,$ otherwise, we get zero.
The second case will give us a thin substitution
to the polynomial $h,$ thus, such summands are also zero. Therefore,
the substitution (\ref{substK22}) to the polynomial $h_{\mu}$ will give
the same result as this substitution to the polynomial
\[h \Bigl(f'_1 \circ f'_2 \circ \zeta_1, f'_3 \circ \zeta_2, \dots,
f'_{p} \circ \zeta_q, \zeta_{q+1},\dots,\zeta_{\tilde{s}}, \
\tilde{x}_1,\dots,\tilde{x}_{n-\tilde{s}} \Bigr),\]
where $f'_l=(\prod_{m=1}^{\tilde{s}+\mu}
\mathcal{A}_{Y_{l, m}} \mathcal{A}_{Z_{l, m}}) f_l.$ Similarly to arguments
of Lemma \ref{ind-simple}, we can see that the result of our substitution
to the polynomial $\tilde{f}'_l=(\prod_{m=1}^{\tilde{s}+\mu}
\mathcal{A}_{Y_{l, m}} \mathcal{A}_{Z_{l, m}}) \tilde{f}_l$ is equal to
$2^{\xi_l} \varepsilon_l.$
Since $\varepsilon_l$ is symmetric, and $f'_l$ is the symmetric
part of $\tilde{f}'_l$ then for $f'_l$ our substitution also yields
$2^{\xi_l} \varepsilon_l,$ and
for the polynomial $h_\mu$ it gives the result
$\alpha h(\varepsilon_{\bf 1} \circ \varepsilon_{\bf 2} \circ
r_{1, ({\bf 1},{\bf 2})}, \varepsilon_{\bf 3} \circ r_{2, ({\bf 3},l_2'')},\dots, \varepsilon_{\bf p}
\circ r_{q, (l'_{q},{\bf p})},\dots, b_1,\dots,b_{n-\tilde{s}})=$
$\alpha h(a_1,\dots,a_n)$ for $\alpha = 2^{\xi},$ \
$\xi=\sum_{l=1}^{p} \xi_l.$
Therefore, $h_\mu$ is the desired polynomial. The substitution $\bar{u}$ is given by
(\ref{substK22}).

The case $p=1$ is
analogous to and simpler than the previous one. If in the case
$p=1$ we have
$(\varepsilon_1/2) \circ h(a_1,\dots,a_n) = h(a_1,\dots,a_n)$ then
we need to consider the substitution (\ref{substK22}) to the
polynomial
\begin{equation*}
h_\mu=
\Bigl( \prod_{m=1}^{\tilde{s}+\mu}
(\mathcal{A}_{\mathcal{Z}^{(0)}_m}  \cdot \mathcal{A}_{\mathcal{Z}^{(1)}_m}) \Bigr)
f_1 \circ
h \Bigl( \zeta_1, \dots,\zeta_{\tilde{s}}, \
\tilde{x}_1,\dots,\tilde{x}_{n-\tilde{s}} \Bigr).
\end{equation*}
Otherwise the substitution $\bar{a}$ contains an element $r_{s, (1,2)}$,
and the construction is similar to the previous case.
In case $p=0$ the algebra $A$ is nilpotent, and we can assume that $h_{\mu}=h,$
and $\bar{u}=\bar{a}.$  \hfill $\Box$

\begin{lemma} \label{gamma}
Let $A$ be a $*$PI-reduced algebra then \
$\mathrm{ind}_*(A)=\mathrm{par}_*(A).$ If $f$ is an exact polynomial
for $A,$ and $f \notin \mathrm{Id}^*(A)$ then $*T[f] \cap
S_{\mu}(A) \ne \emptyset$ for any $\mu \in \mathbb{N}_0.$
\end{lemma}
\noindent {\bf Proof.} By Lemma \ref{Exact1}, $A$ has an exact
polynomial $\tilde{f}$ which is not a $*$-identity of $A.$ Next, $\tilde{f}$ can be
nonzero only for a complete substitution containing exactly
$(\mathrm{nd}(A)-1)$ radical elements. Lemma \ref{Gammasub}
implies that $\tilde{f}$ has a nontrivial consequence
$\widetilde{g} \notin \mathrm{Id}^{*}(A)$ of type
$(\mathrm{dims}_* A;\mathrm{nd}(A)-1;\mu)$ for any $\mu \in
\mathbb{N}_0.$ Since by Lemma \ref{ind-simple} we have
$\beta(A)=\mathrm{dims}_* A$ then Definition \ref{defgamma}
implies that $\gamma(A) > \mathrm{nd}(A)-1.$ Taking into account
$\mathrm{ind}_*(A) \le \mathrm{par}_*(A)$ we obtain
$\gamma(A)=\mathrm{nd}(A).$

Moreover, by Lemma \ref{Gammasub} any exact for $A$ polynomial
$f$ that is not a $*$-identity of $A$ for any $\mu \in
\mathbb{N}_0$ has a nontrivial consequence $g_\mu \in
*T[f],$ where $g_\mu$ is $\mu$-boundary polynomial for $A.$
\hfill $\Box$

Lemma \ref{gamma} together with Lemma \ref{Exact2} immediately
implies

\begin{lemma} \label{Exact3}
Any multilinear $\mu$-boundary polynomial for a $*$PI-reduced
algebra $A$ is exact for $A$.
\end{lemma}

\begin{lemma} \label{Exact4}
Given a $*$PI-reduced algebra $A,$ and a nonnegative integer $\mu,$
let $S_{A,\mu}$ be any set of $*$-polynomials of
type $(\beta(A);\gamma(A)-1;\mu).$ Then any multilinear \mbox{$*$-polynomial}
$f \in *T[S_{A,\mu}]+\mathrm{Id}^{*}(A)$ is exact for $A.$
\end{lemma}

\begin{lemma} \label{Kmu}
Let $\Gamma$ be a proper $*$T-ideal, and $A$ a $*$PI-reduced algebra
such that $\mathrm{ind}_*(\Gamma)=\mathrm{ind}_*(A).$
Suppose that a $*$-polynomial $f$ satisfies the conditions
$f \notin \mathrm{Id}^{*}(A),$ and $f \in
K_{\hat{\mu}}(\Gamma)+\mathrm{Id}^{*}(A)$ for some $\hat{\mu} \in \mathbb{N}_0.$
Then for any $\mu \in \mathbb{N}_0$
we have that $*T[f] \cap S_{\mu}(A) \ne \emptyset.$
\end{lemma}
\noindent {\bf Proof.} The full linearization $\tilde{f} \in F \langle Y, Z \rangle$
of a some multihomogeneous component of $f$ also satisfies
$\tilde{f} \in K_{\hat{\mu}}(\Gamma)+\mathrm{Id}^{*}(A),$ \ $\tilde{f}
\notin \mathrm{Id}^{*}(A).$ Then by Lemma \ref{Exact4},
$\tilde{f}$ is exact for $A.$ Now by Lemma \ref{gamma} we obtain that
$\emptyset \ne *T[\tilde{f}] \cap S_{\mu}(A) \subseteq *T[f]
\cap S_{\mu}(A)$ for any $\mu \in \mathbb{N}_0.$ \hfill $\Box$

\section{Representable algebras.}

Let $R$ be a commutative associative $F$-algebra with unit.
Suppose that an $F$-algebra $A$ with involution has a structure of
$R$-algebra, and the involution of $A$ is $R$-linear, i.e.
$r a=a r,$ \ $(ra)^*=r a^*$ for all $r \in R,$ \ $a \in A.$
Particularly, this happens if $R=F$ or if $R \subseteq Z(A) \cap A^{+},$
where $Z(A)$ is the center of $A,$ and $A^{+}$ is its symmetric part.

\begin{definition}
Any $R$-multilinear mapping $\mathfrak{f}:A^{n} \rightarrow R$ is
called $n$-linear form on an $R$-algebra $A$ with involution.
\end{definition}
We construct an analogue of the free algebra
with forms for algebras with involution (see also \cite{Kem2},
\cite{Rasm1}, \cite{Zubk}).

Let us fix a bilinear form $\mathfrak{f}_2$, and a linear form
$\mathfrak{f}_1.$
Denote by $\mathcal{S}$ the free associative commutative algebra
with unit generated by all symbols $\mathfrak{f}_2(u_1,u_2),$
$\mathfrak{f}_1(u_3),$  where $u_1, u_2, u_3 \in F\langle
X^*\rangle$ are nonempty
associative noncommutative $*$-monomials in $X^{*}$ (or $u_1, u_2, u_3 \in F\langle
Y, Z\rangle$). We say that $FS\langle
X^{*} \rangle = F\langle X^{*}\rangle \otimes_F \mathcal{S}$
($FS\langle Y, Z \rangle = F\langle Y, Z \rangle \otimes_F \mathcal{S}$)
is the free $*$-algebra with forms (the free algebra with forms and involution).
We assume that $(f \otimes s_1)
s_2 = s_2 (f \otimes s_1) = f \otimes (s_1 s_2)$ for all $f \in
F\langle X^{*} \rangle,$ \ $s_1, s_2 \in \mathcal{S}.$ The involution
on $FS\langle X^{*}\rangle$ is induced from $F\langle X^{*}\rangle$
by $(f \otimes s_1)^{*}=f^* \otimes s_1.$ The elements of
$FS\langle X^*\rangle$ are called {\it $*$-polynomials with
forms}. The elements of $\mathcal{S}$ are called {\it pure form
$*$-polynomials}.

Next, $F\langle X^*\rangle \otimes_F 1$ will be identified with
$F\langle X^*\rangle.$ The symbol $\otimes$ will be usually
omitted for $*$-polynomials with forms. The
degrees (homogeneity, multilinearity, alternating etc., respectively) of $*$-polynomials
with forms or pure form \mbox{$*$-polynomials} are defined similarly to the
case of ordinary $*$-polynomials assuming that $\deg_x
\mathfrak{f}_2(u_1,u_2)=\deg_x u_1+\deg_x u_2,$ \  $\deg_x
\mathfrak{f}_1(u_1)=\deg_x u_1$ for any $x \in X$ ($x \in Y \cup Z$).

The $\mathcal{S}$-bilinear form
$\mathfrak{f}_2: FS\langle X^*\rangle^2 \rightarrow \mathcal{S}$
is defined on the $\mathcal{S}$-algebra $FS\langle X^*\rangle$ by
\ $\mathfrak{f}_2(\sum_{i} u_i s_i,\sum_{j} u'_j s'_j)
=\sum_{i, j} \mathfrak{f}_2(u_i,u'_j)
s_i s'_j,$ where $u_i, u'_j \in F\langle X^*\rangle$ are monomials in
$X^{*},$ $s_i, s'_j \in \mathcal{S}.$
The $\mathcal{S}$-linear form
$\mathfrak{f}_1: FS\langle X^*\rangle \rightarrow \mathcal{S}$ is
defined on $FS\langle X^*\rangle$ by $\mathfrak{f}_1(\sum_{i} u_i s_i)
=\sum_{i} \mathfrak{f}_1(u_i) s_i,$ where \  $u_i \in F\langle X^*\rangle$
are monomials in $X^{*},$ $s_i \in \mathcal{S}.$

Let $A$ be an $R$-algebra with involution and forms,
$f(x_{1},\dots,x_{n}) \in FS\langle X^*\rangle$
a $*$-polynomial with forms. Then
$A$ satisfies the $*$-identity with forms $f=0$ if $f(a_{1},\dots,a_{n})=0$
for any $a_{1}, \dots ,a_{n} \in A.$
The ideal of $*$-identities with forms
$\mathrm{SId}^*(A)=\{
f \in FS\langle X^*\rangle | \  A $ satisfies $ f=0 \}$ is
an $\mathcal{S}$-ideal of $FS\langle X^*\rangle$ invariant
with respect to involution and closed
under all $*$-endomorphisms of the algebra $FS\langle X^*\rangle$
which preserve the forms. Observe that $\mathrm{SId}^*(A)$ has the property that
$g_1 \cdot \mathfrak{f}_2(f,g_2),
g_1 \cdot \mathfrak{f}_2(g_2,f), g_1 \cdot \mathfrak{f}_1(f)
\in \mathrm{SId}^*(A)$ \
for any $g_1, g_2 \in FS\langle X^*\rangle,$ \
$f \in \mathrm{SId}^*(A).$
Ideals of $FS\langle X^*\rangle$ with
all mentioned properties are called $*$TS-ideals.
Given a $*$TS-ideal $\widetilde{\Gamma},$ and given polynomials with forms
$f, g \in FS\langle X^*\rangle$ the notation
$f=g \ (\mathrm{mod} \ \widetilde{\Gamma})$ means that $f-g \in
\widetilde{\Gamma}.$ Denote by $*TS[\mathcal{V}]$ the $*$TS-ideal generated by a set
$\mathcal{V} \subseteq FS\langle X^*\rangle.$

Let $\widetilde{\Gamma} \unlhd FS\langle X^*\rangle$ be a $*$TS-ideal of
$FS\langle X^*\rangle$. Denote by $I=\mathrm{Span}_F\{
\mathfrak{f}_2(f,u) v,$ $\mathfrak{f}_2(u,f) v,$  $\mathfrak{f}_1(f) v  |
f \in \widetilde{\Gamma}, u \in F\langle X^*\rangle, v \in \mathcal{S} \}
\unlhd \mathcal{S}$ the ideal of $\mathcal{S}$ generated by
all elements of the forms $\mathfrak{f}_2(f,u),$ $\mathfrak{f}_2(u,f),$
$\mathfrak{f}_1(f),$ \ $f \in \widetilde{\Gamma},$ \  $u \in F\langle X^*\rangle.$
Let $\bar{\mathcal{S}}=\mathcal{S}/I$ be the quotient algebra.
The quotient algebra with involution $\overline{FS}\langle
X^*\rangle=FS\langle X^* \rangle /\widetilde{\Gamma}$ ($\overline{FS}\langle
Y, Z\rangle=FS\langle Y, Z \rangle /\widetilde{\Gamma}$ ) has the
structure of a $\bar{\mathcal{S}}$-algebra. The
$\bar{\mathcal{S}}$-bilinear function $\mathfrak{f}_2:(FS\langle
X^*\rangle / \widetilde{\Gamma})^2 \rightarrow \bar{\mathcal{S}}$
is naturally defined by
$\mathfrak{f}_2(a_1+\widetilde{\Gamma},a_2+\widetilde{\Gamma})=
\mathfrak{f}_2(a_1,a_2) + I,$  \
and the $\bar{\mathcal{S}}$-linear function $\mathfrak{f}_1:(FS\langle
X^*\rangle / \widetilde{\Gamma}) \rightarrow \bar{\mathcal{S}}$
is defined by
$\mathfrak{f}_1(a_1+\widetilde{\Gamma})=
\mathfrak{f}_1(a_1) + I,$  \  $a_1, a_2 \in FS\langle X^*\rangle.$

The ideal of $*$-identities with forms of $\overline{FS}\langle
X^*\rangle$ coincides with $\widetilde{\Gamma}.$ Moreover,
$FS\langle X^*\rangle /\widetilde{\Gamma}$ is the
relatively free $*$-algebra with forms for the
$*$TS-ideal $\widetilde{\Gamma}.$

We also can consider the free associative $*$-algebra with
forms $FS\langle X_{\nu}^*\rangle$ and the relatively free $*$-algebra
with forms $FS\langle X_{\nu}^*\rangle/(\widetilde{\Gamma} \cap FS\langle
X_{\nu}^*\rangle)$ of a finite rank $\nu \in \mathbb{N}.$

Let us take a finite dimensional $F$-algebra with involution $A=B \oplus J$
with the Jacobson radical $J=J(A),$ and the semisimple part $B.$
Consider for any element $b \in B$ the linear operator
$\mathfrak{T}_{b}:B \rightarrow B$  on the $*$-subalgebra $B$
defined by
\begin{eqnarray} \label{Oper}
\mathfrak{T}_b(c)= b \circ c, \quad  c \in B.
\end{eqnarray}
It is clear
that $\mathfrak{T}_{\alpha_1 b_1 + \alpha_2 b_2}=
\alpha_1 \mathfrak{T}_{b_1}+\alpha_2 \mathfrak{T}_{b_2}$ for all
$\alpha_i \in F,$ \ $b_i \in B.$ If $b$ is a symmetric element with respect
to involution then the subspaces $B^{+},$ $B^{-}$
are stable under $\mathfrak{T}_{b}.$
If $b$ is skew-symmetric then $\mathfrak{T}_{b}(B^{+}) \subseteq B^{-},$
and $\mathfrak{T}_{b}(B^{-}) \subseteq B^{+}.$ Particularly, the trace of the
operator $\mathfrak{T}_{b}$ is zero for any skew-symmetric element
$b \in B^{-}.$

Then the bilinear form $\mathfrak{f}_2:A^2 \rightarrow F,$
and the linear form $\mathfrak{f}_1:A \rightarrow F$
are naturally defined on $A$ by the rules
\begin{eqnarray} \label{Atrace}
&&\mathfrak{f}_2(a_1,a_2)=\mathfrak{f}_2(b_1,b_2)=
\mathrm{Tr}(\mathfrak{T}_{b_1} \cdot \mathfrak{T}_{b_2}),\nonumber \\
&&\mathfrak{f}_1(a_1)=\mathfrak{f}_1(b_1)=
\mathrm{Tr}(\mathfrak{T}_{b_1}),
\quad a_i=b_i+r_i \in A, \ b_i \in B, \  r \in J,
\end{eqnarray}
where $\mathfrak{T}_{1} \cdot \mathfrak{T}_{2}$ is the product of linear
operators, and $\mathrm{Tr}$ is the usual trace.
It is clear that $\mathfrak{f}_2$ is symmetric form satisfying
$\mathfrak{f}_2(r,a)=0$ for any $r \in J,$ \  $a \in A,$ and
$\mathfrak{f}_2(a_1,a_2)=0$ for any $a_1 \in A^{-},$ \ $a_2 \in A^{+}.$
The linear form $\mathfrak{f}_1$ also satisfies
$\mathfrak{f}_1(r)=0$ for any $r \in J,$ \
$\mathfrak{f}_1(a)=0$ for any $a \in A^{-}.$
Particularly, the next lemma holds.

\begin{lemma} \label{Traceid10}
A finite dimensional $*$-algebra $A$ with the forms
(\ref{Atrace}) over a field $F$ satisfies $*$-identities with forms \\
$\begin{array}{llll}
\qquad \qquad
&\mathfrak{f}_2(\tilde{y},\tilde{z}) \cdot f=0, \qquad
&\mathfrak{f}_2(\tilde{z},\tilde{y}) \cdot f=0, \qquad
&\mathfrak{f}_1(\tilde{z}) \cdot f=0,
\end{array}$
\\ where $f \in FS\langle Y, Z \rangle$ is any form $*$-polynomial,
\  $\tilde{y} \in Y,$ \  $\tilde{z} \in Z.$
\end{lemma}

\begin{lemma} \label{Traceid1}
Given a finite dimensional $*$-algebra $A$ with the forms
(\ref{Atrace}) over a field $F,$ and a $*$-polynomial
$f \in F\langle Y, Z \rangle$
of type $(\mathrm{dims}_* A,\mathrm{nd}(A)-1,1)$
suppose that $\{y_1,\dots,y_{t_1} \} \subseteq Y$
is the symmetric part, and $\{z_1,\dots,z_{t_2} \} \subseteq Z$
the skew-symmetric part of the set in which $f$ is
$(\mathrm{dims}_* A)$-alternating ($t_1=\dim B^{+},$ \ $t_2=\dim B^{-}$).
Then $A$ satisfies the $*$-identities with forms \\
$\begin{array}{ll}
\qquad
&\mathfrak{f}_2(\tilde{y}_1,\tilde{y}_2)f=
\sum_{i=1}^{t_1} f|_{y_i:=\tilde{y}_1 \circ (\tilde{y}_2 \circ y_i)} +
\sum_{i=1}^{t_2} f|_{z_i:=\tilde{y}_1 \circ (\tilde{y}_2 \circ z_i)}, \quad
\tilde{y}_1, \tilde{y}_2 \in Y, \\
\qquad
&\mathfrak{f}_2(\tilde{z}_1,\tilde{z}_2)f=
\sum_{i=1}^{t_1} f|_{y_i:=\tilde{z}_1 \circ (\tilde{z}_2 \circ y_i)} +
\sum_{i=1}^{t_2} f|_{z_i:=\tilde{z}_1 \circ (\tilde{z}_2 \circ z_i)}, \quad
\tilde{z}_1, \tilde{z}_2 \in Z, \\
\qquad
&\mathfrak{f}_1(\tilde{y}_1)f=
\sum_{i=1}^{t_1} f|_{y_i:=\tilde{y}_1 \circ y_i} +
\sum_{i=1}^{t_2} f|_{z_i:=\tilde{y}_1 \circ z_i}, \quad
\tilde{y} \in Y. \\
\end{array}$
\end{lemma}
\noindent {\bf Proof.} The $*$-polynomials with forms
$g_1=\mathfrak{f}_2(\tilde{y}_1,\tilde{y}_2)f-\sum_{i=1}^{t_1}
f|_{y_i:=\tilde{y}_1 \circ (\tilde{y}_2 \circ y_i)} -
\sum_{i=1}^{t_2} f|_{z_i:=\tilde{y}_1 \circ (\tilde{y}_2 \circ z_i)},$ \
$g_2=\mathfrak{f}_2(\tilde{z}_1,\tilde{z}_2)f-\sum_{i=1}^{t_1} f|_{y_i:=
\tilde{z}_1 \circ (\tilde{z}_2 \circ y_i)} -
\sum_{i=1}^{t_2} f|_{z_i:=\tilde{z}_1 \circ (\tilde{z}_2 \circ z_i)},$
$g_3=\mathfrak{f}_1(\tilde{y}_1)f-\sum_{i=1}^{t_1} f|_{y_i:=\tilde{y}_1
\circ y_i} - \sum_{i=1}^{t_2} f|_{z_i:=\tilde{y}_1 \circ z_i}$ \
are $\tau_j$-alternating in $\mathcal{Z}_j \subseteq Y \cup Z$ with
$\tau_j > \mathrm{dims}_{*} A$ for any $j=1,\dots,\mathrm{nd}(A)-1.$
It is enough to consider an elementary replacement of the variables
$\bigcup_{j=1}^{\mathrm{nd}(A)-1} \mathcal{Z}_j \cup
\{ \tilde{y}_1, \tilde{y}_2, \tilde{z}_1, \tilde{z}_2 \} \cup
\{ y_1,\dots, y_{t_1}, z_1,\dots, z_{t_2}  \}$
of $g_1,$ $g_2,$ $g_3.$
It is clear that at least one variable of any set $\mathcal{Z}_j$ should
be replaced by a radical element, otherwise the result is zero.
Hence, the
variables $y_1, \dots, y_{t_1},$ $z_1,\dots, z_{t_2},$ and
$\tilde{y}_1, \tilde{y}_2$ (or $\tilde{z}_1, \tilde{z}_2$)
must be replaced by semisimple elements only,
and the sets $\{ y_1, \dots, y_{t_1} \},$ $\{ z_1,\dots, z_{t_2} \}$ must be
replaced by pairwise distinct elements of the basis $D_0,$ and $D_1$
(\ref{basisD}) respectively, otherwise $g_i$ also yields zero.

Suppose that a polynomial $\tilde{f}(x_1,\dots,x_{\tilde{t}},\dots)$ is alternating
in a $*$-homogeneous set of variables $x_1,\dots,x_{\tilde{t}} \subset Y$ (or from $Z$),
and $\widetilde{B} \subseteq A$ is a $*$-homogeneous subspace of $A$
of symmetric (or skew-symmetric) elements with
$\dim \widetilde{B} =\tilde{t}.$ It can be directly checked (see also \cite{BelRow},
Theorem J) that for an arbitrary linear operator $\mathfrak{T}:\widetilde{B}
\rightarrow \widetilde{B},$ and for all pairwise distinct basic elements
$b_1,\dots,b_{\tilde{t}}$ of $\widetilde{B}$ the equality
$\mathrm{Tr}(\mathfrak{T}) f(b_1,\dots,b_{\tilde{t}},\dots)=f(\mathfrak{T}(b_1),
\dots,b_{\tilde{t}},\dots)+\cdots+f(b_1,\dots,\mathfrak{T}(b_{\tilde{t}}),\dots)$
holds in $A,$ the replacement of other variables being arbitrary.

The subspaces $B^{+},$ and $B^{-}$ are closed under the actions
of $\mathfrak{T}_{\tilde{b}_1}$ for any $\tilde{b}_1 \in B^{+},$
and
$\mathfrak{T}_{\tilde{b}_1} \mathfrak{T}_{\tilde{b}_2}$ for any
$\tilde{b}_1, \tilde{b}_2 \in B^{+}$ (or $\tilde{b}_1, \tilde{b}_2 \in B^{-}$).
Particularly,
$\mathrm{Tr}(\mathfrak{T}_{\tilde{b}_1})=
\mathrm{Tr}((\mathfrak{T}_{\tilde{b}_1})^{+})+
\mathrm{Tr}((\mathfrak{T}_{\tilde{b}_1})^{-}),$
and
$\mathrm{Tr}(\mathfrak{T}_{\tilde{b}_1} \mathfrak{T}_{\tilde{b}_2})=
\mathrm{Tr}((\mathfrak{T}_{\tilde{b}_1}
\mathfrak{T}_{\tilde{b}_2})^{+})+\mathrm{Tr}((\mathfrak{T}_{\tilde{b}_1}
\mathfrak{T}_{\tilde{b}_2})^{-}),$
where $(\mathfrak{T})^{\delta}$ is the restriction of the operator
$\mathfrak{T}$ on $B^{\delta},$ \  $\delta \in \{ +, - \}.$
The application of the remark above to the linear operators
$(\mathfrak{T}_{\tilde{b}_1})^{+},$ \
$(\mathfrak{T}_{\tilde{b}_1})^{-},$ \
$(\mathfrak{T}_{\tilde{b}_1} \mathfrak{T}_{\tilde{b}_2})^{+},$ and
$(\mathfrak{T}_{\tilde{b}_1} \mathfrak{T}_{\tilde{b}_2})^{-}$
completes the proof.
\hfill $\Box$

\begin{lemma} \label{Traceid2} \qquad
Let $f(\widetilde{x}_1,\dots,\widetilde{x}_k) \in F\langle Y, Z \rangle$
be a $*$-polynomial of a type
$(\beta;\gamma-1;1)$ (for some $\beta \in \mathbb{N}_0^2,$ \ $\gamma \in
\mathbb{N}$), and $s(\zeta_1,\dots,\zeta_d) \in \mathcal{S}$ ($\{\zeta_1,\dots,\zeta_d\} \subseteq
Y \cup Z$) be a pure form $*$-polynomial. Then
there exists a $*$-polynomial
$g_s(\widetilde{x}_1,\dots,\widetilde{x}_k,\zeta_1,\dots,\zeta_d) \in *T[f]$
such that any finite
dimensional $*$-algebra $A$ with forms (\ref{Atrace}) having parameter
$\mathrm{par}_*(A)=(\beta;\gamma)$ satisfies the $*$-identity
with forms
\begin{equation*}
s(\zeta_1,\dots,\zeta_d) \cdot f(\widetilde{x}_1,\dots,\widetilde{x}_k) -
g_s(\widetilde{x}_1,\dots,\widetilde{x}_k,\zeta_1,\dots,\zeta_d) =0.
\end{equation*}
\end{lemma}
\noindent {\bf Proof.} Assume that
$f$ is $(\mathrm{dims}_* A)$-alternating in
$\{ \widetilde{x}_1,\dots,\widetilde{x}_t \},$ \ $t=t_1+t_2.$
Let us take for any $i=1,\dots,t_1$ a right normed jordan monomial $w_i$ of the algebra
$(F\langle Y,Z \rangle^{+},+,\circ)$ of type
$w_i=\zeta'_{i1} \circ (\zeta'_{i 2} \circ (\dots (\zeta'_{i m_{i-1}} \circ \zeta'_{i m_{i}})))$
having even total degree in variables from $Z,$ \
$\zeta'_{i j} \in Y \cup Z,$ \  $m_{i} \in \mathbb{N}.$
Consider also for indices $i=t_1+1,\dots,t_1+t_2$ arbitrary right normed jordan monomials
$w_i \in F\langle Y,Z \rangle^{-}$ in $Y \cup Z$
of odd total degree in $Z.$
Then from Lemma \ref{Traceid1} it follows that
the algebra $A$ satisfies form $*$-identities
\begin{eqnarray*}
\mathfrak{f}_2(\tilde{\zeta}_1,\tilde{\zeta}_2)f(w_1,\dots,w_{t},\widetilde
{X})-(f(\tilde{\zeta_1} \circ (\tilde{\zeta_2} \circ w_1),\dots,w_{t},\widetilde
{X})+\\
\dots+f(w_1,\dots,\tilde{\zeta_1} \circ (\tilde{\zeta_2} \circ  w_{t}),\widetilde{X}))=0, \\
\mathfrak{f}_1(\tilde{y})f(w_1,\dots,w_{t},\widetilde
{X})-(f(\tilde{y} \circ w_1,\dots,w_{t},\widetilde
{X})+\dots+f(w_1,\dots,\tilde{y} \circ w_{t},\widetilde{X}))=0,
\end{eqnarray*}
where  $\tilde{\zeta}_1, \tilde{\zeta}_2 \in Y$
or $\tilde{\zeta}_1, \tilde{\zeta}_2 \in Z,$ and
$\tilde{y} \in Y.$

By induction on the total number $n=n_1+n_2$ of forms
$\mathfrak{f}_i$ we obtain the next equality
$\mathfrak{f}_2(\tilde{\zeta}_{1},\tilde{\zeta}_{2}) \cdots
\mathfrak{f}_2(\tilde{\zeta}_{2 n_2-1},\tilde{\zeta}_{2 n_2})
\mathfrak{f}_1(\tilde{\zeta}_{2 n_2+1}) \cdots
\mathfrak{f}_1(\tilde{\zeta}_{2 n_2+n_1})
f(w_1,\dots,w_{t},\widetilde{X})=\sum_{l=1}^{\widetilde{n}}
f(\widetilde{w}_{l 1},\dots,\widetilde{w}_{l t},\widetilde{X})$
$(\mathrm{mod} \ \mathrm{SId}^*(A)),$ where $w_i,$ \ $\widetilde{w}_{l i}$
are right normed jordan symmetric or skew-symmetric monomials, and $\tilde{\zeta}_{j} \in Y \cup Z$
are the corresponding symmetric or skew-symmetric variables.
Hence, the equality
$\mathfrak{f}_2(u_1,u_2) \cdots \mathfrak{f}_2(u_{2 n_2-1},u_{2 n_2}) \times$
$\mathfrak{f}_1(u_{2 n_2+1}) \cdots \mathfrak{f}_1(u_{2 n_2+n_1})
f(\widetilde{x}_1,\dots,\widetilde{x}_{k})=
g_{(u)}(\widetilde{x}_1,\dots,\widetilde{x}_{k},\zeta_1,\dots,\zeta_d)$
$(\mathrm{mod} \ \mathrm{SId}^*(A))$ holds for all $*$-homogeneous elements
$u_{2i-1}, u_{2i} \in F\langle Y,Z \rangle^{\delta},$ \ $\delta \in \{ +, - \}$
($i=1,\dots,n_2$),
$u_{i} \in F\langle Y,Z \rangle^{+}$ ($i=2 n_2+1,\dots,2 n_2+n_1$),
where $g_{(u)} \in *T[f]$ is some $*$-polynomial without forms.
Observe that the structure
of $g_{(u)}$ depends only on the polynomials $u_i,$
and $\zeta_1,\dots,\zeta_d$ are
the variables of the polynomial $u_1 \cdots u_{2n_2+n_1}.$

Any pure form polynomial with involution $s \in \mathcal{S}$ can be written in the form
$s(\zeta_1,\dots,\zeta_d) \ = \ \sum_{(j)} \ \alpha_{(j)} \
\mathfrak{f}_2(u_{j_1},u_{j_2}) \cdots
\mathfrak{f}_2(u_{j_{2n_2-1}},u_{j_{2n_2}})
\mathfrak{f}_1(u_{j_{2n_2+1}}) \cdots
\mathfrak{f}_1(u_{j_{2n_2+n_1}}),$ \\ where
$u_{j_l}$ are $*$-homogeneous elements of $F\langle Y,Z \rangle$
(symmetric or skew-symmetric), and $\alpha_{(j)} \in F.$
Taking into account Lemma \ref{Traceid10} and arguments above we obtain that
$A$ satisfies the form $*$-identity
$s(\zeta_1,\dots,\zeta_d) \cdot
f(\widetilde{x}_1,\dots,\widetilde{x}_{k}) - g_s(\widetilde{x}_1,\dots,
\widetilde{x}_{k},\zeta_1,\dots,\zeta_d)=0,$ where
$g_s \in *T[f]$ is a $*$-polynomial which does not depend on $A.$
\hfill $\Box$

Let $A$ be a finite dimensional $*$-algebra over a field $F$
satisfying the claims of Lemma \ref{Pierce}.
Assume that $B$ is its semisimple part, and $J$ is the Jacobson radical.
Then denote $\dim B^{+} = t_1,$ $\dim B^{-} = t_2,$ $\dim J^{+} = q_1,$
$\dim J^{-} = q_2.$
We take the set
$\Lambda_{\nu}=\{ \lambda_{\theta i j} | \theta =1, 2; 1 \le i \le
\nu; 1 \le j \le t_\theta+q_\theta \}$ for any positive integer $\nu$.
Consider the free commutative associative algebra
with unit $F[\Lambda_{\nu}]^{\#}$ generated by
$\Lambda_{\nu}$, and the associative algebra
$\mathcal{P}_{\nu}(A)=F[\Lambda_{\nu}]^{\#} \otimes_F A.$
The algebra $\mathcal{P}_{\nu}(A)$ has the structure of a $F[\Lambda_{\nu}]
^{\#}$-module defined by \ $a \cdot f = f \cdot a  = f \cdot (\sum_{i} f_i
\otimes a_i)=\sum_{i} (f f_i) \otimes a_i,$ \ for any $f, f_i \in
F[\Lambda_{\nu}]^{\#},$ \ $a_i \in A,$ \ $a=\sum_{i} f_i \otimes
a_i \in \mathcal{P}_{\nu}(A).$ The involution on $\mathcal{P}_{\nu}(A)$
is naturally induced from $A$ by $(f \otimes a)^*=f \otimes a^*,$ \
$f \in F[\Lambda_{\nu}]^{\#},$ \ $a \in A.$

We define the $F[\Lambda_{\nu}]^{\#}$-bilinear map
$\mathfrak{f}_2: \mathcal{P}_{\nu}(A)^2 \rightarrow
F[\Lambda_{\nu}]^{\#},$ and the $F[\Lambda_{\nu}]^{\#}$-linear map
$\mathfrak{f}_1: \mathcal{P}_{\nu}(A) \rightarrow
F[\Lambda_{\nu}]^{\#}$
\begin{eqnarray} \label{Ptrace}
&&\mathfrak{f}_2(a,a')=\mathfrak{f}_2(\sum_{i} f_{i} \otimes a_{i},
\sum_{j} f'_j \otimes a'_{j})=\sum_{i,j} f_i f'_j\
\mathfrak{f}_2(a_i,a'_j), \nonumber \\
&&\mathfrak{f}_1(a)=\mathfrak{f}_1(\sum_{i} f_{i} \otimes a_{i})=
\sum_{i} f_i \ \mathfrak{f}_1(a_i).
\end{eqnarray}
They are well defined bilinear and linear forms on $\mathcal{P}_{\nu}(A),$
respectively. Here
$a_{i}, a'_j \in A,$ \
$f_i, f'_j \in F[\Lambda_{\nu}]^{\#},$ and $\mathfrak{f}_2(a_i,a'_j),$
$\mathfrak{f}_1(a_i)$ are the forms (\ref{Atrace}) defined on $A.$

We call by {\it a Cayley-Hamilton type
$*$-\!polynomial} a degree homogeneous $*$-po\-ly\-no\-mial with forms of the following type
\begin{eqnarray*}
x^{n}+\sum \limits_{\mathop{\ i_0 + i_1 + \dots + j_{k_2} =n,}
\limits_{\scriptstyle 0 < i_0 < n, \  1 \le k_2+k_1}}
\alpha_{(i),(j)} \ \ x^{i_0} \mathfrak{f}_2(x^{i_1},x^{j_1}) \cdots
\mathfrak{f}_2(x^{i_{k_2}},x^{j_{k_2}})
\mathfrak{f}_1(x^{i_{k_2+1}}) \cdots \mathfrak{f}_1(x^{i_{k_2+k_1}}),
\end{eqnarray*}
where $\alpha_{(i),(j)} \in F.$  Note that here $i_l, j_l > 0$ ($l \ge 0$).

For example, it is well known (\cite{Proc1}, \cite{Rasm1})
that the full matrix algebra $M_n(F)$
satisfies the Cayley-Hamilton identity with trace $\mathcal{X}_{n}(x)=0,$
which is defined recurrently by
the formulas $\mathcal{X}_{1}(x)=x-\mathrm{Tr}(x),$ \quad
$\mathcal{X}_{n}(x)=\mathcal{X}_{n-1}(x)\cdot x - \frac{1}{n} \cdot
\mathrm{Tr}(\mathcal{X}_{n-1}(x)\cdot x),$ where
$\mathrm{Tr}$ is the usual trace of matrix. Hence, $M_n(F)$
with any involution and with forms
$\mathfrak{f}_2(a,b)=\mathrm{Tr}(a \cdot b),$ \
$\mathfrak{f}_1(a)=\mathrm{Tr}(a)$
($a, b \in M_n(F)$) satisfies the Cayley-Hamilton type $*$-identity with
forms $\widetilde{\mathcal{X}}_{n}(x) \cdot x=0,$ where
$\widetilde{\mathcal{X}}_{n}(x)$ is the $*$-polynomial with forms
obtained from $\mathcal{X}_{n}(x)$ by exchange
$\mathrm{Tr}(x^s)=\mathfrak{f}_1(x^s),$ \ $s \ge 1.$

\begin{lemma} \label{HamKel}
$\mathcal{P}_{\nu}(A)$ satisfies a Cayley-Hamilton type $*$-identity
$\mathcal{K}_{3t+1}^{\mathrm{nd}(A)}(x)=0$ for some
Cayley-Hamilton type $*$-polynomial $\mathcal{K}_{3t+1}(x)$
of degree $3t+1,$ \ $t=t_1+t_2.$
\end{lemma}
\noindent {\bf Proof.}
Let us consider the semisimple part $B$ of $A$ with operation $\circ.$
Then $B$ is a finite dimensional Jordan algebra, $\dim B=t=t_1+t_2.$
It is well known that this Jordan algebra satisfies
the next relation
\begin{eqnarray} \label{pol-Op}
\mathfrak{T}_{b^n} = \sum_{0 \le k \le n/2} \  \alpha_k \  \mathfrak{T}^k_{b^2}
\cdot \mathfrak{T}^{n-2k}_{b}, \quad  b \in B, \  \  n \in \mathbb{N},
\end{eqnarray}
where the operators were defined by (\ref{Oper}), numbers $\alpha_k \in \mathbb{Q}$
are fixed, that do not depend on $b.$ It is clear that both of the operators $\mathfrak{T}_{b^2},$ and
$\mathfrak{T}_{b}$ satisfy the Cayley-Hamilton identity with trace of degree
$t.$ If $n=3 t$ then either $k \ge t$ or $n-2k \ge t.$ Therefore,
either $\mathfrak{T}^k_{b^2}$ or $\mathfrak{T}^{n-2k}_{b}$ is expressed
via the powers of lesser degrees and their traces. Thus, by applying the operator
$\mathfrak{T}_{b^{3t}}$ to the element $b,$ we obtain
\begin{eqnarray} \label{pol-Op1}
b^{3t+1} = \sum
\limits_{\mathop{ m + 2(i_1 + \dots + i_{k'}) + }
\limits_{\mathop{ j_1 + \dots + j_{k''}=3t+1, }
\limits_{\scriptstyle 1 \le m < 3t+1, \  1 \le k'+k'' }}}
\alpha'_{m, (i),(j)} \
b^{m} \  \mathrm{Tr}( \mathfrak{T}^{i_1}_{b^2}) \cdots
\mathrm{Tr}(\mathfrak{T}^{i_{k'}}_{b^2}) \mathrm{Tr}(\mathfrak{T}^{j_1}_{b})
\cdots \mathrm{Tr}(\mathfrak{T}^{j_{k''}}_{b})
\end{eqnarray}
for any $b \in B,$ where the coefficients $\alpha'_{m,(i),(j)} \in \mathbb{Q}$
do not depend on $b.$ It can be directly checked also that
for any $s \ge 1$ we have
\[ \mathrm{Tr}(\mathfrak{T}^{s}_{c})=\alpha''_{s,0} \ \
\mathfrak{f}_1(c^s) \  +
\sum_{1 \le l \le s/2}  \alpha''_{s,l} \  \
\mathfrak{f}_2(c^{l},c^{s-l})
\]
for some rational coefficients $\alpha''_{s,l} \in \mathbb{Q},$
that do not depend on $c,$ and for any $c \in B.$
Particularly, this is true for
$c=b,$ and $c=b^2,$ \ $b \in B.$ Thus, it follows from
(\ref{pol-Op1}) that the associative $*$-algebra $B$ with forms
(\ref{Atrace}) satisfies the Cayley-Hamilton type $*$-identity with forms
$\mathcal{K}_{3t+1}(x)=0$
of degree $3t+1.$

Since $A=B \oplus J,$ $J$ is nilpotent of class $\mathrm{nd}(A),$
and the forms $\mathfrak{f}_2,$ \ $\mathfrak{f}_1$ can be non-zero only on semisimple elements,
we conclude that $A$ satisfies $\mathcal{K}_{3t+1}^{\mathrm{nd}(A)}(x)=0.$
Next, $F[\Lambda_{\nu}]^{\#}$ is commutative non-nilpotent algebra,
the involution of $\mathcal{P}_{\nu}(A)$ is $F[\Lambda_{\nu}]^{\#}$-linear,
and the forms (\ref{Ptrace}) are $F[\Lambda_{\nu}]^{\#}$-multilinear.
Therefore it is clear that $\mathcal{P}_{\nu}(A)=F[\Lambda_{\nu}]^{\#} \otimes_F A$
also satisfies $\mathcal{K}_{3t+1}^{\mathrm{nd}(A)}(x)=0.$
\hfill $\Box$

Let  $\{\hat{b}_{1 1}, \dots,
\hat{b}_{1 t_1}\}$ be a basis of the symmetric part
$B^{+}$ of the semisimple part $B$ of $A,$ \ $\dim B^{+} =
t_1,$ and $\{\hat{b}_{2 1}, \dots,
\hat{b}_{2 t_2}\}$ a basis of the skew-symmetric part
$B^{-}$ of $B,$  \ $\dim B^{-} =t_2,$ \ $t=t_1+t_2.$
Also let
$\{\hat{r}_{1 1}, \dots, \hat{r}_{1
q_1}\}$ be a basis of the symmetric part $J^{+}$ of the
Jacobson radical $J=J(A),$  \ $\dim J^{+} = q_1,$
and $\{\hat{r}_{2 1}, \dots, \hat{r}_{2
q_2}\}$ a basis of the skew-symmetric part $J^{-}$ of $J(A),$ correspondingly,
\ $\dim J^{-} = q_2.$
Recall that all these bases can be chosen belonging to the set $D_0 \cup D_1 \cup U_0 \cup U_1$
((\ref{basisD}), (\ref{basisU}), Lemma \ref{Pierce}). Take
the elements
\begin{eqnarray} \label{genset}
\mathfrak{y}_{\theta i}=\sum_{j=1}^{t_\theta} \lambda_{\theta i j} \otimes
\hat{b}_{\theta j} + \sum_{j=1}^{q_\theta} \lambda_{\theta i
j+t_\theta} \otimes \hat{r}_{\theta j} \  \in  \mathcal{P}_{\nu}(A),
\qquad \theta=1,2, \  \ 1 \le i \le \nu.
\end{eqnarray}
All $\mathfrak{y}_{\theta i}$ are $*$-homogeneous (symmetric for $\theta=1$,
and skew-symmetric for $\theta=2$). Consider for any positive integer $\nu$
the $F$-subalgebra $\mathcal{F}_{\nu}(A)=\langle \mathfrak{y}_{\theta i} |
\theta=1,2; \ 1 \le i \le \nu \rangle$ of $\mathcal{P}_{\nu}(A)$
generated by  $\mathfrak{Y}^{*}_{\nu}=\{ \mathfrak{y}_{\theta i} |
\theta=1,2; \  1 \le i \le \nu \}.$ Then $\mathcal{F}_{\nu}(A)$ is
$*$-invariant. Any map $\varphi$ of the generators to
arbitrary $*$-homogeneous elements $\widetilde{a}_{\theta i} \in A$ \
($\widetilde{a}_{1 i} \in A^{+},$ \
$\widetilde{a}_{2 i} \in A^{-},$ \
$\widetilde{\alpha}_{\theta i j} \in F$)
\begin{equation} \label{hom1}
\varphi:\mathfrak{y}_{\theta i} \mapsto \widetilde{a}_{\theta
i}=\sum_{j=1}^{t_\theta} \widetilde{\alpha}_{\theta i j}
\hat{b}_{\theta j} + \sum_{j=1}^{q_\theta}
\widetilde{\alpha}_{\theta i j+t_\theta} \hat{r}_{\theta j} \quad (\theta=1,2; \quad i=1,\dots,\nu)
\end{equation}
can be extended to the $*$-homomorphism of $F$-algebras
$\varphi:\mathcal{F}_{\nu}(A) \rightarrow A$, also inducing the
$*$-homomorphism $\widetilde{\varphi}:\mathcal{P}_{\nu}(A)
\rightarrow A$ defined by the following equalities
\begin{eqnarray} \label{hom2}
\widetilde{\varphi}((\lambda_{\theta_1 i_1 j_1} \cdots
\lambda_{\theta_k i_k j_k})\otimes a)=(\widetilde{\alpha}_{\theta_1
i_1 j_1} \cdots \widetilde{\alpha}_{\theta_k i_k j_k}) \cdot a
\qquad \forall a \in A.
\end{eqnarray}
Observe that the $*$-homomorphism $\widetilde{\varphi}$ preserves the
forms defined by (\ref{Ptrace}) on $\mathcal{P}_{\nu}(A)$ and
by (\ref{Atrace}) on $A.$

The elements of $\mathcal{F}_{\nu}(A)$
are called {\it quasi-polynomials in the variables $\mathfrak{Y}^{*}_{\nu}.$} Products
of the generators $\mathfrak{y}_{\theta i} \in \mathfrak{Y}^{*}_{\nu}$ of the algebra
$\mathcal{F}_{\nu}(A)$ are called {\it quasi-monomials}. We have
also that $\mathrm{Id}^*(\mathcal{F}_{\nu}(A)) \supseteq
\mathrm{Id}^*(\mathcal{P}_{\nu}(A))=\mathrm{Id}^{*}(A)$ for any $\nu
\in \mathbb{N}.$

Recall that $\mathcal{F}_{\nu}(A)$ is a finitely generated PI-algebra.
By Shirshov's height theorem   \cite{Shirsh} $\mathcal{F}_{\nu}(A)$
has a finite height and a finite Shirshov's basis.
More precisely, there exist an integer
$\mathcal{H},$ and elements $w_1,\dots,w_d \in \mathcal{F}_{\nu}(A)$ such that
any element $u \in \mathcal{F}_{\nu}(A)$ has the form
$u=\sum_{(i)=(i_1,\dots,i_k)}
 \alpha_{(i)} \  w_{i_1}^{c_1} \dots w_{i_k}^{c_k},$
where $k \le \mathcal{H},$ \  $\{ i_1,\dots,i_k \} \subseteq
\{1,\dots,d\},$ \  $c_j \in \mathbb{N},$ \  $\alpha_{(i)} \in F.$

Consider the polynomials
$\hat{\mathfrak{s}}_{i,(l_1,l_2)}=
\mathfrak{f}_2(w_i^{l_1},w_i^{l_2}) \in
F[\Lambda_{\nu}]^{\#}$ ($i=1,\dots,d,$ \
$l_1,l_2=1,\dots,3t$),
and
$\hat{\mathfrak{s}}_{i,l}=
\mathfrak{f}_1(w_i^{l}) \in
F[\Lambda_{\nu}]^{\#}$
($i=1,\dots,d,$ \ $l=1,\dots,3t$).
Then
$\widehat{F}=F[
\hat{\mathfrak{s}}_{i,(l_1,l_2)},
\hat{\mathfrak{s}}_{i l_1} \ | \ 1 \le i \le d; \ 1
\le l_1, l_2 \le 3t \ ]^{\#}$ is
the associative commutative $F$-subalgebra
of $F[\Lambda_{\nu}]^{\#}$ with unit
generated by $\{ \hat{\mathfrak{s}}_{i,(l_1,l_2)},
\hat{\mathfrak{s}}_{i l} \}$, and by
the unit of $F[\Lambda_{\nu}]^{\#}.$

Take the $*$-invariant $\widehat{F}$-subalgebra
$\mathcal{T}_{\nu}(A)=\widehat{F} \mathcal{F}_{\nu}(A)$ of
$\mathcal{P}_{\nu}(A).$ Then $\mathcal{F}_{\nu}(A)$ is a
$*$-subalgebra of $\mathcal{T}_{\nu}(A).$ An arbitrary
map of type (\ref{hom1}) can be uniquely extended to a
$*$-homomorphism from $\mathcal{T}_{\nu}(A)$ to $A$ preserving the
forms (it is the restriction of
$\widetilde{\varphi}$ defined by (\ref{hom2}) onto $\mathcal{T}_{\nu}(A)$).
It follows from Lemma \ref{HamKel} that all elements $w_i$
are algebraic of degree $\mathrm{nd}(A)(3t+1)$ over
$\widehat{F}.$
Therefore, by the Shirshov's height theorem,
$\mathcal{T}_{\nu}(A)$ is a finitely generated
$\widehat{F}$-module, where $\widehat{F}$ is Noetherian. By
theorem of Beidar \cite{Beid} the algebra $\mathcal{T}_{\nu}(A)$
is representable.

Let $V \subseteq F\langle Y, Z \rangle$ be a set of $*$-polynomials.
We denote by $V(\mathcal{T}_{\nu}(A)) \unlhd
\mathcal{T}_{\nu}(A)$ the verbal $*$-ideal generated by results of all
$*$-homogeneous substitutions of elements of $\mathcal{T}_{\nu}(A)$
to any $*$-polynomial from $V.$

\begin{remark} \label{VerbId}
Given a set $V \subseteq F\langle Y, Z \rangle,$ and a
positive integer $\nu,$ the verbal $*$-ideal $V(\mathcal{T}_{\nu}(A))$ is
$\widehat{F}$-closed. The quotient algebra
$\overline{\mathcal{T}}_{\nu}(A,V)=\mathcal{T}_{\nu}(A)/V(\mathcal{T}_{\nu}(A))$
is representable $\widehat{F}$-algebra with involution.
The ideal of $*$-identities of $\overline{\mathcal{T}}_{\nu}(A,V)$ satisfies \\
\ $\mathrm{Id}^*(\overline{\mathcal{T}}_{\nu}(A,V)) \supseteq
\mathrm{Id}^*(A)+*T[V]+\mathrm{Id}^* \Bigl(F\langle Y_{\nu}, Z_{\nu}
\rangle/\bigl((\mathrm{Id}^*(A)+*T[V]) \cap F\langle Y_{\nu}, Z_{\nu}
\rangle \bigr) \Bigr).$
\end{remark}

\begin{definition} \label{Homog}
We say that a subset $V \subseteq F\langle Y, Z \rangle$ is
$*$-multihomogeneous if $V$ is $*$-invariant, and for any $f \in V$
it contains all multihomogeneous components of $f.$
\end{definition}

\begin{lemma} \label{GS1}
Let $A=A_1 \times \dots \times A_\rho$
be the direct product of arbitrary finite dimensional $*$-algebras
$A_1, \dots, A_\rho.$ Suppose that
$V \subseteq F\langle Y, Z \rangle$ is a $*$-multihomogeneous set, and
$\nu >0$ is an integer. Let us take any $f(\zeta_1,\dots,\zeta_n)
\in \mathrm{Id}^*(\overline{\mathcal{T}}_{\nu}(A,V))$ \
($\zeta_i \in Y \cup Z$), and any
$*$-homogeneous $*$-polynomials $h_1,\dots,h_n \in F\langle
Y_{\nu}, Z_{\nu} \rangle$ ($h_l$ symmetric or skew-symmetric
in accordance with $\zeta_l,$ \ $l=1,\dots,n$).
Then the equality
$f(h_1, \dots,h_n)=\sum_{j} \mathfrak{s}_j
\cdot v_{j 1} \tilde{f}_j(u_{j 1},\dots,u_{j m}) v_{j 2} \ (\mathrm{mod } \
\mathrm{SId}^{*}(A_i))$ holds for
any $i=1,\dots,\rho.$ Here $\tilde{f}_j$ are the full
linearizations of some polynomials $f_j \in V;$ \ $u_{j l} \in
F\langle Y_{\nu}, Z_{\nu} \rangle$ are $*$-homogeneous $*$-polynomials;
$v_{j l} \in F\langle Y_{\nu}, Z_{\nu} \rangle$ are
monomials, possibly empty; and $\mathfrak{s}_j
\in \mathcal{S}$ are pure form $*$-polynomials in the
variables $Y_{\nu} \cup Z_{\nu}.$
\end{lemma}
\noindent {\bf Proof.}
Given a $*$-polynomial $f(\zeta_1,\dots,\zeta_n) \in
\mathrm{Id}^*(\overline{\mathcal{T}}_{\nu}(A,V)),$ and arbitrary
$*$-ho\-mo\-ge\-neous polynomials $h_1,\dots,h_n \in F\langle Y_{\nu}, Z_{\nu}
\rangle,$ symmetric or skew-symmetric according to variables of $f,$
we have $f(\tilde{h}_1,\dots,\tilde{h}_n) \in V(\mathcal{T}_{\nu}(A)).$ Here
the quasi-polynomial $\tilde{h}_i=h_i(\mathfrak{y}_{1 1},\dots,\mathfrak{y}_{2 \nu})$
is obtained by the replacement of the variables $Y_{\nu} \cup Z_{\nu}$
by the corresponding elements of $\mathfrak{Y}^{*}_{\nu}.$
Hence, in the algebra $\mathcal{T}_{\nu}(A)$ we obtain the equality
$f(\tilde{h}_1,\dots,\tilde{h}_n)=
\sum_{j} \widetilde{\mathfrak{s}}_j \cdot \widetilde{v}_{j 1}
\tilde{f}_j(\widetilde{u}_{j 1},\dots,\widetilde{u}_{j n}) \widetilde{v}_{j 2},$ where
$\tilde{f}_j$ are the full linearizations of polynomials $f_j
\in V,$ \ $\widetilde{u}_{j l}=u_{j l}(\mathfrak{y}_{1 1},\dots,\mathfrak{y}_{2 \nu})
\in \mathcal{F}_{\nu}(A)$ are $*$-homogeneous quasi-polynomials, \
$\widetilde{\mathfrak{s}}_j=\mathfrak{s}_j(\mathfrak{y}_{1 1},\dots,\mathfrak{y}_{2 \nu})
\in \widehat{F}$ are pure form quasi-polynomials,
$\widetilde{v}_{j l}=v_{j l}(\mathfrak{y}_{1 1},\dots,\mathfrak{y}_{2 \nu}) \in
\mathcal{F}_{\nu}(A)$ are quasi-monomials, possibly empty.

An arbitrary map $\varphi:\mathfrak{Y}^{*}_{\nu} \rightarrow
A_i$ from the generating set (\ref{genset}) of
$\mathcal{F}_{\nu}(A)$ into any subalgebra $A_i \subseteq A$
($i=1,\dots,\rho$) can be extended to the $*$-homomorphism
$\widetilde{\varphi}:\mathcal{T}_{\nu}(A) \rightarrow A_i$ preserving
forms. Then an equality of quasi-polynomials in the algebra
$\mathcal{T}_{\nu}(A,V)$ implies the equivalence of the
corresponding form $*$-polynomials with respect to the
$*$TS-ideal $\mathrm{SId}^{*}(A_i).$ \hfill $\Box$

\begin{lemma} \label{GS}
Suppose that $A_1,\dots,A_\rho$ are
any $*$PI-reduced algebras  with
$\mathrm{ind}_*(A_i)=\kappa$ for all $i=1,\dots,\rho.$ Given a
subset $V \subseteq \bigcup_{i=1}^{\rho} S_{\mu}(A_i)$ (for any
$\mu \ge 1$), and a positive integer $\nu,$ there exists a $F$-finite
dimensional $*$-algebra $C_\nu$ such
that $\mathrm{Id}^{*}(C_\nu) = \mathrm{Id}^{*}\Bigl(F\langle
Y_{\nu}, Z_{\nu} \rangle/\bigl((\cap_{i=1}^\rho
\mathrm{Id}^{*}(A_i)+*T[V]) \cap F\langle Y_{\nu}, Z_{\nu} \rangle
\bigr) \Bigr).$
\end{lemma}
\noindent {\bf Proof.}
Let us take $A=A_1 \times \dots \times A_\rho.$
By Lemmas \ref{GS1}, \ref{Traceid2}, for any $f(\zeta_1,\dots,\zeta_n) \in
\mathrm{Id}^*(\overline{\mathcal{T}}_{\nu}(A,V)),$ and for any
$*$-homogeneous $*$-polynomials $h_1,\dots,h_n \in F\langle Y_{\nu}, Z_{\nu}
\rangle$ (symmetric or skew-symmetric according to variables of $f,$ \
$i=1,\dots,n$) we have
\begin{eqnarray*}
&&f(h_1, \dots,h_n)=\sum_{j} \mathfrak{s}_j v_{j 1}
\tilde{f}_j(u_{j 1},\dots,u_{j n}) v_{j 2}
(\mathrm{mod } \ \mathrm{SId}^{*}(A_i))=\\
&&\sum_{j} v_{j 1} \tilde{g}_j v_{j 2} \ (\mathrm{mod } \
\mathrm{SId}^{*}(A_i)) \qquad \mbox{ for any } i=1,\dots,\rho.
\end{eqnarray*}
Here $\tilde{f}_j
\in *T[V];$ $\tilde{g}_j \in *T[\tilde{f}_j] \cap F\langle
Y_{\nu}, Z_{\nu} \rangle;$ $u_{j l} \in F\langle
Y_{\nu}, Z_{\nu} \rangle$ are $*$-homogeneous $*$-polynomials, $v_{j l} \in
F\langle Y_{\nu}, Z_{\nu} \rangle$ are monomials, possibly empty;
$\mathfrak{s}_j \in \mathcal{S}.$  Then $g=\sum_{j}
v_{j 1} \tilde{g}_j v_{j 2} \in *T[V]
\cap F\langle Y_{\nu}, Z_{\nu} \rangle.$ Therefore $f(h_1,\dots,h_n)-g \in$
$\mathrm{SId}^*(A_i) \cap F\langle Y_{\nu}, Z_{\nu} \rangle =$
$\mathrm{Id}^{*}(A_i) \cap F\langle Y_{\nu}, Z_{\nu} \rangle$ for all
$i=1,\dots,\rho.$ Hence $f(h_1,\dots,h_n) \in \bigl(\mathrm{Id}^{*}(A)
+ *T[V] \bigr) \cap F\langle Y_{\nu}, Z_{\nu} \rangle,$ \ and $f \in
\mathrm{Id}^{*} \Bigl(F\langle Y_{\nu}, Z_{\nu} \rangle/\bigl((\mathrm{Id}^{*}(A)+*T[V])
\cap F\langle Y_{\nu}, Z_{\nu}
\rangle \bigr) \Bigr).$

By Remark \ref{VerbId} we have also that \quad $\mathrm{Id}^{*} \Bigl(
F\langle Y_{\nu}, Z_{\nu}
\rangle/\bigl((\mathrm{Id}^{*}(A)+*T[V]) \cap F\langle Y_{\nu}, Z_{\nu}
\rangle \bigr) \Bigr) \subseteq$
$\mathrm{Id}^{*}(\overline{\mathcal{T}}_{\nu}(A,V)).$ The algebra
$\overline{\mathcal{T}}_{\nu}(A,V)$ is representable.
Hence from Lemma \ref{Repr} it follows that there exists a finite dimensional
over $F$ $*$-algebra $C_\nu$ such that
$\mathrm{Id}^{*}(C_\nu)=\mathrm{Id}^{*}(\overline{\mathcal{T}}_{\nu}(A,V))=\mathrm{Id}^{*}
\Bigl(F\langle Y_{\nu}, Z_{\nu} \rangle/\bigl((\mathrm{Id}^{*}(A)+*T[V])
\cap F\langle Y_{\nu}, Z_{\nu} \rangle \bigr) \Bigr).$ \hfill $\Box$

\section{Identities with involution of finitely generated algebras.}

\begin{lemma} \label{Main}
Let $\Gamma$  be a non-trivial ideal of $*$-identities of
a finitely generated associative
$*$-algebra over a field $F$ of zero characteristic.
Then there exists a finite dimensional associative
$F$-algebra with involution $\widetilde{A}$ satisfying the conditions
$\mathrm{Id}^{*}(\widetilde{A}) \subseteq
\Gamma,$ \ $\mathrm{ind}_*(\widetilde{A})=\mathrm{ind}_*(\Gamma),$ and
$S_{\hat{\mu}}(\mathcal{O}(\widetilde{A})) \cap \Gamma =\emptyset$ for some
$\hat{\mu} \in \mathbb{N}_0.$
\end{lemma}
\noindent {\bf Proof.} By Lemma \ref{lemma1}, $\Gamma$
contains the ideal of $*$-identities of some
finite dimensional $*$-algebra $A.$ By Lemma
\ref{Smu}, we can suppose that $A=\mathcal{O}(A) \times \mathcal{Y}(A)$ with the
senior components $A_1,\dots,A_{\rho}.$ It is clear that
$\kappa=\mathrm{ind}_*(\Gamma) \le \kappa_1=\mathrm{ind}_*(A)$
(Lemma \ref{Prop1}). If $\Gamma \subseteq \mathrm{Id}^{*}(A_{i})$
for some $i=1,\dots,\rho$ then $\kappa_1 =
\mathrm{ind}_*(A_{i})=\kappa.$ Thus, the case $\Gamma \subseteq
\mathrm{Id}^{*}(\mathcal{O}(A))$ is trivial.

Assume that $\Gamma \nsubseteq \mathrm{Id}^{*}(A_i)$ for all
$i=1,\dots,\rho_1$ ($1 \le \rho_1 \le \rho$), and $\Gamma \subseteq
\mathrm{Id}^{*}(A''),$ where $A''=\times_{i=\rho_1+1}^\rho A_i,$ \
$A'=\times_{i=1}^{\rho_1} A_i.$ Consider the set
$V=S_{\tilde{\mu}}(A') \cap \Gamma$ for
$\tilde{\mu}=\widehat{\mu}(A)$ (Definition
\ref{sen}). Take $\nu=\mathrm{rk}(D)$ for a finitely generated
$*$-algebra $D$ such that $\Gamma=\mathrm{Id}^{*}(D).$ By
Lemma \ref{GS}, there exists a finite dimensional $F$-algebra
with involution $C_\nu$ such that
$\mathrm{Id}^{*}(C_\nu)=\mathrm{Id}^{*} \bigl(F\langle X_{\nu}^{*}
\rangle/\bigl((\mathrm{Id}^{*}(A')+*T[V]) \bigcap F\langle
X_{\nu}^{*} \rangle \bigr) \bigr).$

Then $\widetilde{A}=C_\nu \times A'' \times \mathcal{Y}(A)$ is
a finite dimensional algebra with involution. For any $f(x_1,\dots,x_n)
\in \mathrm{Id}^{*}(\widetilde{A}),$ and for all $*$-polynomials
$h_1,\dots,h_n \in F\langle X_\nu^{*} \rangle$
we have $f(h_1,\dots,h_n)=h+g,$ where $h
\in \mathrm{Id}^{*}(A'),$ \ $g \in \Gamma \cap
\mathrm{Id}^{*}(\mathcal{Y}(A)).$ Therefore, $h=f(h_1,\dots,h_n)-g \in
\mathrm{Id}^{*}(A') \cap \mathrm{Id}^{*}(A'') \cap
\mathrm{Id}^{*}(\mathcal{Y}(A))=\mathrm{Id}^{*}(A) \subseteq \Gamma.$
Hence, $f(h_1,\dots,h_n)=h+g \in \Gamma$ for any $h_1,\dots,h_n \in
F\langle X_{\nu}^{*} \rangle,$ and
$\mathrm{Id}^{*}(\widetilde{A}) \subseteq \Gamma$ by virtue of Remark \ref{freefg}.
Particularly,
$\mathrm{ind}_*(\widetilde{A}) \ge \kappa.$

Suppose that $A'' \ne 0.$ Then $\kappa_1=
\kappa=\mathrm{ind}_{*}(\widetilde{A}).$ Thus, either $\mathrm{ind}_*(C_\nu)=\kappa$
implying $\mathcal{O}(\widetilde{A})=A'' \times \mathcal{O}(C_\nu),$
or $\mathrm{ind}_*(C_\nu) < \kappa$ implying
$\mathcal{O}(\widetilde{A})=A''.$  Since, $S_\mu(A'') \cap \Gamma=\emptyset$ then in the
first case we use Lemmas \ref{Prop1}, \ref{subset}, \ref{Prop2}, \ref{Smu} and conclude that
\ for any $\mu \ge \max \{ \widehat{\mu}(C_\nu), \tilde{\mu} \}$ we have
$S_\mu(\mathcal{O}(\widetilde{A})) \cap \Gamma=$ $S_\mu(C_\nu)
\cap \Gamma \subseteq$ $S_\mu(A') \cap \Gamma \subseteq$ $V
\subseteq \mathrm{Id}^{*}(C_\nu).$ Thus, in both cases
$S_\mu(\mathcal{O}(\widetilde{A})) \bigcap \Gamma = \emptyset,$
and $\widetilde{A}$ satisfies the claims of the Lemma.

If $A'' = 0$ then $\kappa \le \mathrm{ind}_{*}(\widetilde{A})=$
$\max \{ \mathrm{ind}_*(C_\nu), \mathrm{ind}_*(\mathcal{Y}(A)) \}
\le$ $\kappa_1.$  If $\mathrm{ind}_*(C_\nu) =
\kappa_1=\mathrm{ind}_*(A')$ then $\mathcal{O}(\widetilde{A})=
\mathcal{O}(C_\nu),$ and by analogy with previous case we apply
Lemmas \ref{subset}, \ref{Smu} for any $\mu \ge \max \{
\widehat{\mu}(C_\nu), \tilde{\mu} \}$ and obtain
$S_\mu(\mathcal{O}(\widetilde{A})) \cap \Gamma = \emptyset.$
It also follows from Lemma \ref{subset}  that
$\mathrm{ind}_*(\widetilde{A}) = \mathrm{ind}_*(\Gamma),$ and
$\widetilde{A}$ is also the desired algebra.

The last case $A'' = 0,$ $\mathrm{ind}_*(C_\nu) < \kappa_1$ gives
$\widetilde{A}=C_\nu \times \mathcal{Y}(A)$ with
$\mathrm{Id}^{*}(\widetilde{A}) \subseteq \Gamma,$ and $\kappa \le
\mathrm{ind}_*(\widetilde{A}) < \kappa_1=\mathrm{ind}_*(A).$ Then by
the inductive step on $\mathrm{ind}_*(A)$ we can assume that the
assertion of the Lemma holds in this case. \hfill $\Box$

Lemma \ref{Main}, and Lemma \ref{subset} imply the
corollary.
\begin{lemma} \label{SmuG}
Let $\Gamma$ be the non-trivial ideal of $*$-identities of a finitely generated associative
algebra with involution over a field of characteristic zero.
Then $\Gamma$ contains the ideal of $*$-identities of some finite dimensional
associative $*$-algebra $A$ satisfying
$\mathrm{ind}_*(\Gamma)=\mathrm{ind}_*(A),$ and $S_{\tilde{\mu}}(A) =
S_{\tilde{\mu}}(\Gamma)$ for some $\tilde{\mu} \in \mathbb{N}_0$.
\end{lemma}

\begin{theorem} \label{*PI}
A non-zero $*$T-ideal of identities with involution of a finitely generated
associative algebra with involution over a field $F$ of
characteristic zero coincides
with the $*$T-ideal of identities with involution of some $F$-finite dimensional
associative algebra with involution.
\end{theorem}
\noindent {\bf Proof.} Let $\Gamma \ne (0)$ be the ideal of $*$-identities of a
finitely generated algebra with involution $D$. We use the induction on the
Kemer index $\mathrm{ind}_*(\Gamma)=\kappa=(\beta;\gamma)$ of $\Gamma.$

\underline{The base of induction.} \  Let
$\mathrm{ind}_*(\Gamma)=(\beta;\gamma)$ with $\beta=(0,0).$
Then $D$ is a nilpotent finitely generated algebra. Hence, $D$ is finite dimensional.

\underline{The inductive step.}
Lemmas \ref{Main}, \ref{SmuG} imply
$\Gamma \supseteq \mathrm{Id}^{*}(A),$ where
$A=\mathcal{O}(A)+\mathcal{Y}(A)$ is a finite dimensional $*$-algebra
with $\mathrm{ind}_*(\Gamma)=\mathrm{ind}_*(A)=\kappa.$ Moreover,
$S_{\tilde{\mu}}(\Gamma)=S_{\tilde{\mu}}(\mathcal{O}(A))=S_{\tilde{\mu}}(A)
\subseteq \mathrm{Id}^{*}(\mathcal{Y}(A))$ for some $\tilde{\mu} \in
\mathbb{N}_0.$

Let $A_1,\dots,A_{\rho}$ be the senior components of the algebra
$A.$ We can assume that all $A_i$ satisfy the claims of Lemma \ref{Pierce}.
Denote  $(t_1,t_2)=\beta(\Gamma)=\mathrm{dims}_* A_i,$
\ $t=t_1+t_2;$ \ $\gamma=\gamma(\Gamma)=\mathrm{nd}(A_i)$
(for all $i=1,\dots,\rho$). Let us take for any $i=1,\dots,\rho$
the algebra $\widetilde{A}_i=\mathcal{R}_{q_i,s}(A_i)$ defined by
(\ref{FRad}) for the senior component $A_i$ with $q_i=\dim_F A_i,$
\  $s=(t+1)(\gamma+\tilde{\mu}).$ Then $\widetilde{A}_i$ is
a finite dimensional $*$-algebra. We have also
$\Gamma_i=\mathrm{Id}^{*}(\widetilde{A}_i)=$ $\mathrm{Id}^{*}(A_i),$ and
$\mathrm{dims}_* \widetilde{A}_i=$ $\mathrm{dims}_* A_i=$ $\beta.$
The Jacobson radical $J(\widetilde{A}_i)=(X_{q_i}^{*})/I$ is
nilpotent of class at most $s=(t+1)(\gamma+\tilde{\mu}),$
where $I=\Gamma_i(B_i(X_{q_i}^{*}))+(X_{q_i}^{*})^s.$ Here the algebra
$B_i$ can be considered as the semisimple parts of $A_i$ and of
$\widetilde{A}_i$ simultaneously (Lemma \ref{Aqs}, Lemma
\ref{gamma}). Particularly, the algebras $\widetilde{A}_i$ also satisfy
the claims of Lemma \ref{Pierce}.
By Remark \ref{VerbId}, and Lemma \ref{Repr} there
exists a finite dimensional over $F$ algebra with involution $C$ such that
$\mathrm{Id}^{*}(C)=\mathrm{Id}^{*}(\overline{\mathcal{T}}_{\nu}
(\widetilde{A},\Gamma)),$ where $\widetilde{A}=\times_{i=1}^\rho
\widetilde{A}_i,$ $\nu=2 \mathrm{rkh}(D).$

Let us denote $\widetilde{D}_{\nu}=F\langle X_{\nu}^{*}
\rangle/\bigl((\Gamma + K_{\widetilde{\mu}}(\Gamma)) \cap
F\langle X_{\nu}^{*} \rangle \bigr).$ Lemmas \ref{Prop1},
\ref{Prop2} imply that
$\mathrm{ind}_*(\widetilde{D}_{\nu}) \le \mathrm{ind}_*(\Gamma +
K_{\tilde{\mu}}(\Gamma)) < \mathrm{ind}_*(\Gamma)$. By the inductive
step we obtain
$\mathrm{Id}^{*}(\widetilde{D}_{\nu})=\mathrm{Id}^{*}(\widetilde{U}),$
where $\widetilde{U}$ is a finite dimensional over $F$ $*$-algebra.
Remark \ref{VerbId} yields $\Gamma
\subseteq \mathrm{Id}^{*}(C \times \widetilde{U}).$

Consider a multilinear polynomial $f(\tilde{x}_1,\dots,\tilde{x}_m) \in
\mathrm{Id}^{*}(C \times \widetilde{U})$ in symmetric and skew-symmetric
variables $\tilde{x}_i \in Y \cup Z.$
Let us take any multihomogeneous (with respect to degrees of variables)
and $*$-homogeneous $*$-polynomials
$w_1, \dots, w_m \in F\langle Y_{\nu}, Z_{\nu} \rangle$
($w_i$ is symmetric or skew-symmetric according to $\tilde{x}_i,$ \ $i=1,\dots,m$).
We have $f(w_1,\dots,w_m)=g+h$ for some multihomogeneous
$*$-polynomials $g \in \Gamma,$ \ $h \in K_{\tilde{\mu}}(\Gamma)$ also depending
on $Y_{\nu} \cup Z_{\nu}.$
Then by Lemma \ref{GS1} we obtain $h=f(w_1,\dots,w_m)-g \in \mathcal{S}
\Gamma + \mathrm{SId}^{*}(\widetilde{A}_i)$ for any
$i=1,\dots,\rho.$ Hence
$\tilde{h}(\tilde{x}_1,\dots,\tilde{x}_n) \in \mathcal{S} \Gamma +
\mathrm{SId}^{*}(\widetilde{A}_i)$ holds also for the full
linearization $\tilde{h}$ of $h$.

Suppose that
$\bar{a}=(r_1,\dots,r_{\tilde{s}},b_{\tilde{s}+1},\dots,b_n)$ is
any elementary complete substitution of elements of the algebra
$A_i$ in $\tilde{h}$, where $r_j \in J(A_i),$ $b_j \in B_i,$
$\tilde{s}=\gamma-1.$
By Lemmas \ref{Exact4}, \ref{Gammasub} there exists a polynomial
$h_{\tilde{\mu}}(\mathcal{Z}_1,\dots,\mathcal{Z}_{\tilde{s}+\tilde{\mu}},
\mathcal{X}) \in \mathcal{S} \Gamma +
\mathrm{SId}^{*}(\widetilde{A}_i)$ of the respective type, and
an elementary substitution $\bar{u}$ from $A_i$ in $h_{\tilde{\mu}}$
such that $h_{\tilde{\mu}}(\bar{u})=\alpha \tilde{h}(\bar{a}),$
$\alpha \in F,$ \ $\alpha \ne 0.$
Moreover $h_{\tilde{\mu}}$ is alternating in any
$\mathcal{Z}_j$ ($j=1,\dots,\tilde{s}+\tilde{\mu}$), and
all variables from $\mathcal{X}$
are replaced by semisimple elements.
Then we have
\begin{eqnarray} \label{lhdec}
\alpha_2 h_{\tilde{\mu}}=\bigl(
\prod_{m=1}^{\tilde{s}+\tilde{\mu}}
(\mathcal{A}_{\mathcal{Z}^{(0)}_m}  \cdot \mathcal{A}_{\mathcal{Z}^{(1)}_m})
 \bigr) h_{\tilde{\mu}}=
\sum_j  \bigl( \prod_{m=1}^{\tilde{s}+\tilde{\mu}}
(\mathcal{A}_{\mathcal{Z}^{(0)}_m}  \cdot \mathcal{A}_{\mathcal{Z}^{(1)}_m})
\bigr) \bigl(
\mathfrak{\tilde{s}}_j \tilde{g}_j \bigr) (\mathrm{mod } \
\mathrm{SId}^{*}(\widetilde{A}_i)),
\end{eqnarray}
where $\mathcal{Z}^{(0)}_m$ is the subset of symmetric variables of $\mathcal{Z}_m,$ \
$\mathcal{Z}^{(1)}_m$ the subset of skew-symmetric variables; \
$\alpha_2 \in F,$ $\alpha_2 \ne 0;$ \ $\tilde{g}_j
\in \Gamma,$ \ $\mathfrak{\tilde{s}}_j \in \mathcal{S}.$ Denote by
$\{\zeta_1,\dots,\zeta_{\hat{n}} \}$
the variables $\mathcal{Z} \cup \mathcal{X}$ of $h_{\tilde{\mu}}$  (the first variables are
from $\mathcal{Z} = \bigcup_{m=1}^{\tilde{s}+\tilde{\mu}}
\mathcal{Z}_m$).

Let us take in the algebra $\widetilde{A}_i$ elements
$\bar{x}_{\pi(k)}=x_{\pi(k)}+I,$ where
$x_{\pi(k)} \in X_{q_i}$ are variables, $\pi(k)$ is the number of the basic elements
$u_k=c_{\pi(k)}$ in the basis $\{c_1,\dots,c_{q_i} \}$ of $A_i$ that consists of
$*$-homogeneous elements chosen in (\ref{basisD}),
(\ref{basisU}), $1 \le \pi(k) \le q_i.$

Consider the following substitution from $\widetilde{A}_i$ in
$h_{\tilde{\mu}}(\zeta_1,\dots,\zeta_{\hat{n}}).$
\begin{eqnarray} \label{Ai}
&& \zeta_k=(\varepsilon_{l'} \bar{x}_{\pi(k)}
\varepsilon_{l''}+\varepsilon_{l''} \bar{x}^{*}_{\pi(k)}
\varepsilon_{l'})/2 \in J(\widetilde{A}_i) \quad \mbox{ if } \ \ \zeta_k
\in \mathcal{Z}, \  \zeta_k \ \mbox{ is symmetric, and }\nonumber \\
&& \qquad \qquad \qquad \qquad \qquad \qquad  \qquad \qquad \qquad  \quad
u_k=c_{\pi(k)} \in \varepsilon_{l'} A_i
\varepsilon_{l''}+\varepsilon_{l''} A_i
\varepsilon_{l'}; \nonumber \\
&& \zeta_k=(\varepsilon_{l'} \bar{x}_{\pi(k)}
\varepsilon_{l''}-\varepsilon_{l''} \bar{x}^{*}_{\pi(k)}
\varepsilon_{l'})/2 \in J(\widetilde{A}_i) \quad \mbox{ if } \ \  \zeta_k
\in \mathcal{Z}, \  \zeta_k \ \mbox{ is skew-symmetric, }\nonumber \\
&& \qquad \qquad \qquad \qquad \qquad \qquad  \qquad \qquad \qquad  \quad
u_k=c_{\pi(k)} \in \varepsilon_{l'} A_i
\varepsilon_{l''}+\varepsilon_{l''} A_i
\varepsilon_{l'}; \nonumber \\
&& \zeta_k=u_k
\quad \mbox{ if } \ \zeta_k \in \mathcal{X}.
\end{eqnarray}

Suppose that in (\ref{lhdec}) the pure form polynomial $\mathfrak{\tilde{s}}_j$
depends essentially on $\mathcal{Z}.$ Then we get
$\mathfrak{\tilde{s}}_j|_{(\ref{Ai})}=0,$ because the forms on radical
elements are zero (\ref{Atrace}). If $\mathfrak{\tilde{s}}_j$ does not depend
on $\mathcal{Z}$ then $\bigl(
\prod_{m=1}^{\tilde{s}+\tilde{\mu}}
(\mathcal{A}_{\mathcal{Z}^{(0)}_m}  \cdot \mathcal{A}_{\mathcal{Z}^{(1)}_m})
\bigr) \bigl(
\mathfrak{\tilde{s}}_j \tilde{g}_j \bigr) =$
$\mathfrak{\tilde{s}}_j \tilde{\tilde{g}}_j,$ where
$\tilde{\tilde{g}}_j=\bigl(
\prod_{m=1}^{\tilde{s}+\tilde{\mu}}
(\mathcal{A}_{\mathcal{Z}^{(0)}_m}  \cdot \mathcal{A}_{\mathcal{Z}^{(1)}_m})
\bigr)
\tilde{g}_j \in \Gamma.$
If $\tilde{\tilde{g}}_j|_{(\ref{Ai})} \ne 0$ in $\widetilde{A}_i$ then one
of degree multihomogeneous components of
$\tilde{\tilde{g}}_j$ is a $\tilde{\mu}$-boundary polynomial for
$\widetilde{A}_i.$ This implies that $S_{\tilde{\mu}}(A) \cap \Gamma \ne
\emptyset,$ which contradicts to the properties of $A.$ Therefore,
$\tilde{\tilde{g}}_j|_{(\ref{Ai})}=0.$ Thus, in any case
$h_{\tilde{\mu}}|_{(\ref{Ai})}=0$ holds in the algebra
$\widetilde{A}_i.$ Hence, the
substitution
\begin{eqnarray*}
&&\mbox{ if } \ \ \zeta_k
\in \mathcal{Z}, \  \zeta_k \ \mbox{ is symmetric, \  and  }
u_k=c_{\pi(k)} \in \varepsilon_{l'} A_i
\varepsilon_{l''}+\varepsilon_{l''} A_i
\varepsilon_{l'} \quad \mbox{then }  \nonumber \\
&& \zeta_k=v_k=(\varepsilon_{l'} x_{\pi(k)}
\varepsilon_{l''}+\varepsilon_{l''} x^{*}_{\pi(k)}
\varepsilon_{l'})/2; \nonumber \\
&&\mbox{ if } \ \ \zeta_k
\in \mathcal{Z}, \  \zeta_k \ \mbox{ is skew-symmetric, \  and  }
u_k=c_{\pi(k)} \in \varepsilon_{l'} A_i
\varepsilon_{l''}+\varepsilon_{l''} A_i
\varepsilon_{l'} \quad \mbox{then }  \nonumber \\
&& \zeta_k=v_k=(\varepsilon_{l'} x_{\pi(k)}
\varepsilon_{l''}-\varepsilon_{l''} x^{*}_{\pi(k)}
\varepsilon_{l'})/2; \nonumber \\
&& \zeta_k=v_k=u_k
\quad \mbox{ if } \ \zeta_k \in \mathcal{X}
\end{eqnarray*}
in the polynomial $h_{\tilde{\mu}}$ gives the
result $h_{\tilde{\mu}}(v_1,\dots,v_{\hat{n}}) \in
I=\Gamma_i(B_i(X_{q_i}^{*}))+(X_{q_i}^{*})^s$ in the algebra
$B_i(X_{q_i}^{*}).$ Since $|\mathcal{Z}| < s$ then we obtain
$h_{\tilde{\mu}}(v_1,\dots,v_{\hat{n}}) \in
\Gamma_i(B_i(X_{q_i}^{*})).$

Consider the map $\varphi:x_{\pi(k)} \mapsto
c_{\pi(k)}$ if $c_{\pi(k)} \in D_0 \cup D_1$ is a semisimple element, and
$\varphi:x_{\pi(k)} \mapsto
r$ if $c_{\pi(k)}= (\varepsilon_{l'} r \varepsilon_{l''} + (-1)^{\delta}
\varepsilon_{l''} r^{*} \varepsilon_{l'})/2 \in U_0 \cup U_1$ is a radical element
($k=1,\dots,|\mathcal{Z}|$). It is clear that $\varphi$ can be extended to a
$*$-homomorphism $\varphi:B_i(X_{q_i}^{*}) \rightarrow A_i$
assuming $\varphi(b)=b$ for any $b \in B_i.$ Then
$\varphi(h_{\tilde{\mu}}(v_1,\dots,v_{\hat{n}}))=$
$h_{\tilde{\mu}}(\varphi(v_1),\dots,\varphi(v_{\hat{n}}))=$
$h_{\tilde{\mu}}(\bar{u})=$ $\alpha \tilde{h}(\bar{a})
\in $ $\varphi(\Gamma_i(B_i(X_{q_i}^{*})))=(0).$

Therefore $\tilde{h}(\bar{a})=0$ holds in $A_i$ for any elementary
complete substitution $\bar{a} \in A_i^n$ containing $\gamma-1$
radical elements. Since $\tilde{h}$ is
exact for $A_i$ (Lemma \ref{Exact4}), and $\gamma=\mathrm{nd}(A_i)$
then $\tilde{h} \in
\mathrm{Id}^{*}(A_i).$ Hence $h \in \cap_{i=1}^\rho
\mathrm{Id}^{*}(A_i),$ and $h \in
\mathrm{Id}^{*}(\mathcal{O}(A) \times
\mathcal{Y}(A))=\mathrm{Id}^{*}(A) \subseteq \Gamma.$ Thus
we have $f(w_1,\dots,w_m)=g+h \in \Gamma$ for any multihomogeneous
$*$-polynomials $w_1, \dots, w_m \in F\langle Y_{\nu}, Z_{\nu} \rangle,$
which are symmetric or skew-symmetric respectively.
The application of Remark \ref{freefg} now
implies that $\mathrm{Id}^{*}(C \times \widetilde{U}) \subseteq \Gamma.$

Therefore, $\Gamma=\mathrm{Id}^{*}(C \times \widetilde{U}).$
Theorem is proved. \hfill $\Box$

The author is grateful to A.Giambruno and S.P.Mishchenko for inspiration.

\bibliographystyle{amsplain}

\end{document}